\documentclass[a4paper,fleqn]{cas-sc}

\usepackage{setspace}
\usepackage{bm}
\usepackage{multirow}
\usepackage{booktabs}
\usepackage{tabularx}
\usepackage{array}
\usepackage[numbers,sort&compress]{natbib}
\usepackage{graphicx}
\usepackage{float}
\usepackage{geometry}
\usepackage{enumitem}
\usepackage{amsmath}
\usepackage{siunitx}
\usepackage{times}
\usepackage{caption}
\usepackage{subcaption}
\usepackage{rotating}
\usepackage[ruled]{algorithm2e}
\usepackage[title]{appendix}%
\usepackage{placeins}


\usepackage{lineno}
\def\tsc#1{\csdef{#1}{\textsc{\lowercase{#1}}\xspace}}
\tsc{WGM}
\tsc{QE}
\tsc{EP}
\tsc{PMS}
\tsc{BEC}
\tsc{DE}
\begin{document}

\let\WriteBookmarks\relax
\def\floatpagepagefraction{1}
\def\textpagefraction{.001}
\shorttitle{LLE-Guided Midcourse Correction of Low-Energy Earth-Moon Transfers}
\shortauthors{Zihang Wang et~al.}

\title [mode = title]{Local Lyapunov Exponent-Guided Stochastic Midcourse Correction Design of Low-Energy Earth-Moon Transfers}

\tnotemark[1]
\tnotetext[1]{This paper is an extended version of
	Paper IAC-26-A3,IP,195,x112266, presented at the 77th International
	Astronautical Congress, Antalya, T\"urkiye, 5--9 October 2026.}
	
\author[1]{Zihang Wang}[orcid=0009-0007-0345-3063,style=chinese]
\ead{zihangw@buaa.edu.cn}
\credit{Data curation, Formal analysis, Methodology, Software, Writing - Original draft preparation, Writing - review $\&$ editing}
\affiliation[1]{organization={School of Astronautics, Beihang University},
                addressline={Xueyuan Road No.37}, 
                postcode={100191}, 
                city={Beijing},
                country={People's Republic of China}}
\affiliation[2]{organization={State Key Laboratory of High-Efficiency Reusable Aerospace Transportation Technology}, 
	city={Beijing},
	country={People's Republic of China}}

\author[1]{Shuyue Fu}[orcid=0009-0001-5111-9779,style=chinese]
\ead{fushuyue@buaa.edu.cn}
\credit{Conceptualization, Methodology, Writing - review $\&$ editing}

\author[1,2]{Di Wu}[style=chinese]
\ead{wudi2025@buaa.edu.cn}
\credit{Conceptualization, Funding acquisition, Writing - review $\&$ editing}

\author[1,2]{Shengping Gong}[style=chinese]
\cormark[1]
\ead{gongsp@buaa.edu.cn}
\credit{Conceptualization, Funding acquisition, Writing - review $\&$ editing}

\cortext[cor1]{Corresponding author}

\begin{abstract}
Midcourse correction is essential for reliable cislunar transfer missions. For low-energy Earth-Moon transfers, strongly nonuniform error growth makes conventional correction strategies less effective. This paper proposes a stochastic midcourse correction design for low-energy Earth-Moon transfers based on local Lyapunov exponent (LLE). The local error-amplification characteristics of the nominal trajectory are first analyzed using an instantaneous LLE. The analysis shows that strong sensitivity is concentrated near periapsis, and that the LLE maxima coincide with regions of strong finite-time error growth. These maxima are therefore used to segment the trajectory and define the number of corrections. Based on this segmentation, a reduced parameterization is introduced to describe the correction epochs, while second-order state transition tensors (STTs) are employed for nonlinear correction targeting. The resulting timing variables are then optimized through a sample-average approximation of the chance-constrained stochastic problem. Paired Monte Carlo simulations show that, under the same mission-accuracy requirements, the proposed strategy greatly improves the reliability and robustness of trajectory recovery while requiring less correction cost than the classical Breakwell spacing-ratio strategy.
\end{abstract}

\begin{highlights}
\item The instantaneous local Lyapunov exponent (LLE) locates error-amplification regions.
\item LLE-guided segmentation defines correction architecture and reduces timing variables.
\item Chance-constrained scheduling enables reliable, fuel-efficient trajectory recovery.
\end{highlights}


\maketitle
\section{Introduction}\label{sec1}
Cislunar space is becoming an increasingly important domain for future space exploration and operations. As cislunar activities become more frequent, efficient and reliable Earth-Moon transportation have attracted sustained interest~\cite{song2026survey}. Low-energy transfers provide an attractive option by exploiting the multi-body dynamics of cislunar system. These trajectories can reduce propellant demand and offer greater flexibility in flight duration and departure windows~\cite{parker2013targeting,oshima2019long}. In the Earth-Moon circular restricted three-body problem (CRTBP), invariant manifolds associated with libration-point dynamics provide natural pathways for low-energy transport~\cite{koon2001low,howell2006transfers,mingotti2012distant}. Other studies have explored resonant transfers, temporary capture, and numerically optimized trajectories, revealing a wide variety of low-energy transfer solutions in three-body dynamics~\cite{fantino2010note,vaquero2014resonant,dealmeida2026,wang2025}. Four-body models that account for the gravitational influence of the Sun have further been used to  investigate weak stability boundary (WSB) transfers and lunar ballistic capture~\cite{fu2026wsb, martinelli2025wsb,fu2026energy}, expanding the available pathways for low-energy Earth-Moon transportation.  Such multi-body transfer concepts have been applied to lunar mission design, as illustrated by Hiten's ballistic lunar capture~\cite{belbruno1993sun,uesugi1991hiten} and the low-energy transfer developed for EQUULEUS~\cite{oguri2020equuleus}. However, trajectories designed in multi-body dynamical models can exhibit strong sensitivity to small perturbations. Large local gravity gradients near Earth perigees and lunar perilunes can cause rapid amplification of trajectory errors. Even small injection, navigation and maneuver-execution errors may accumulate and lead to substantial deviations from the nominal trajectory \cite{vetrisano2012navigating,qi2018correction,li2025nonsmooth}. Consequently, midcourse correction (MCC) is essential for suppressing accumulated deviations and maintaining the nominal transfer trajectory. Effective MCC is key to the success of cislunar missions. 

The design of MCC involves two closely related aspects: guidance algorithm and correction strategy. The guidance algorithm determines the velocity increment required at each correction epoch. The correction strategy determines how many corrections should be performed and when they should occur along the transfer. For correction strategy design, Breakwell~\cite{breakwell1962spacing} derived the classical spacing-ratio principle for interplanetary navigation. Successive correction epochs are arranged according to their remaining time to the terminal target, providing a simple rule for distributing corrective maneuvers over the flight. Ni et al.~\cite{ni2016breakwell} subsequently examined the effects of the spacing ratio and maneuver number on
correction capability and terminal accuracy. The spacing-ratio principle is also applied to low-energy Earth-Moon transfers by He and Xu~\cite{he2007lowenergy}. To account for various uncertainties, recent studies have also formulated midcourse correction strategy design as a stochastic optimization problem. Geller et al.~\cite{geller2022robust} jointly considered midcourse correction and optical-navigation scheduling in robust cislunar trajectory design. Related studies have further demonstrated the potential of stochastic optimization in spacecraft trajectory design considering uncertainty~\cite{benedikter2022covariance,zhang2023chance,kumagai2025covariance}. For guidance algorithm, higher-order state transition tensors have been introduced specifically in low-energy transfer trajectories to improve nonlinear targeting accuracy. Boone and McMahon ~\cite{boone2021orbital} introduced higher-order state transition tensors (STTs) for nonlinear orbital guidance.  Yin et al.~\cite{yin2023midcourse}  applied high-order STTs to midcourse correction guidance of Earth-Moon distant retrograde orbit transfers. Despite these advances, systematic MCC strategies tailored to the dynamical characteristics of low-energy Earth-Moon transfers remain limited. Existing guidance methods mainly address maneuver computation for prescribed correction epochs and targets. For strategy design, global timing rules do not explicitly account for local sensitivity, while stochastic optimization often assumes a predefined maneuver sequence. As a result, corrections may not be scheduled around critical sensitive passages, allowing trajectory errors to grow and increasing subsequent recovery demands. A systematic approach is therefore needed to incorporate local dynamical sensitivity into the design of correction opportunities and intermediate targets.

Accordingly, this paper proposes a segmented mid-course correction framework for low-energy Earth-Moon transfers based on the local Lyapunov exponent (LLE). The central idea is to construct the correction strategy from the nonuniform dynamical sensitivity of the nominal trajectory and then optimize the remaining timing variables under uncertainty. In this framework, local sensitivity analysis defines where corrections are needed, stochastic optimization determines when they are performed, and nonlinear targeting computes the required correction impulses. The main contributions are summarized as follows. \textbf{First, an instantaneous LLE-based criterion is established to identify local error-amplification regions along low-energy Earth-Moon transfers.}
The LLE is derived from the PCRTBP variational dynamics and analytically related to the local gravity gradient. Local-restart Monte Carlo simulations further verify that the LLE maxima correspond to regions of strong finite-time error amplification, providing a point-wise indicator for locating dynamically sensitive passages. \textbf{Second, a dynamics-informed reduction of the mixed-integer stochastic correction problem is developed.} Jointly designing the maneuver count and correction epochs involves both integer and continuous decision variables. The LLE-based segmentation, together with a prescribed maneuver-placement rule, determines the maneuver count and intermediate targets before optimization. This removes the need to search over the integer variable. A segment-wise timing parameterization further reduces the continuous schedule to one design variable per segment. \textbf{Third, a chance-constrained scheduling procedure is developed to optimize correction timing under uncertainty.} The reduced timing variables are optimized to reduce correction cost subject to a prescribed joint mission-success probability under injection, navigation, and maneuver-execution uncertainties. A sample-average approximation is employed for optimization, and the retained schedules are assessed using independent reliability-validation samples. The joint mission-success criterion incorporates full-path recovery alongside terminal accuracy and maneuver feasibility. Paired Monte Carlo comparisons with the conventional Breakwell strategy are conducted on six representative transfers using the same targeting law and uncertainty model. The proposed method increases the joint mission-success rate by 11.8--90.9 percentage points. Among realizations successfully corrected by both methods, the mean total correction $\Delta v$ is reduced by 11.1--34.2\%.

The remainder of this paper is organized as follows. Section~\ref{sec2} defines the PCRTBP dynamical model, reference low-energy transfers, stochastic correction problem, and performance measures. Section~\ref{sec3} derives the instantaneous LLE, analyzes the local error-amplification mechanism, and develops the sensitivity-guided trajectory segmentation. Section~\ref{sec4} constructs the segmented correction architecture and presents the reduced stochastic schedule optimization and second-order targeting procedure. Section~\ref{sec5} evaluates the correction mechanism and compares the proposed strategy with the Breakwell method through paired Monte Carlo simulations. Section~\ref{sec6} summarizes the main conclusions.

\section{Dynamical Model and Correction Problem}\label{sec2}
The planar Earth--Moon dynamics, reference transfers, uncertainty model, and performance measures are defined below. 

\subsection{Planar Circular Restricted Three-Body Dynamics}\label{subsec2.1}

\begin{figure}[htbp]
	\centering
	\includegraphics[width=0.55\textwidth]{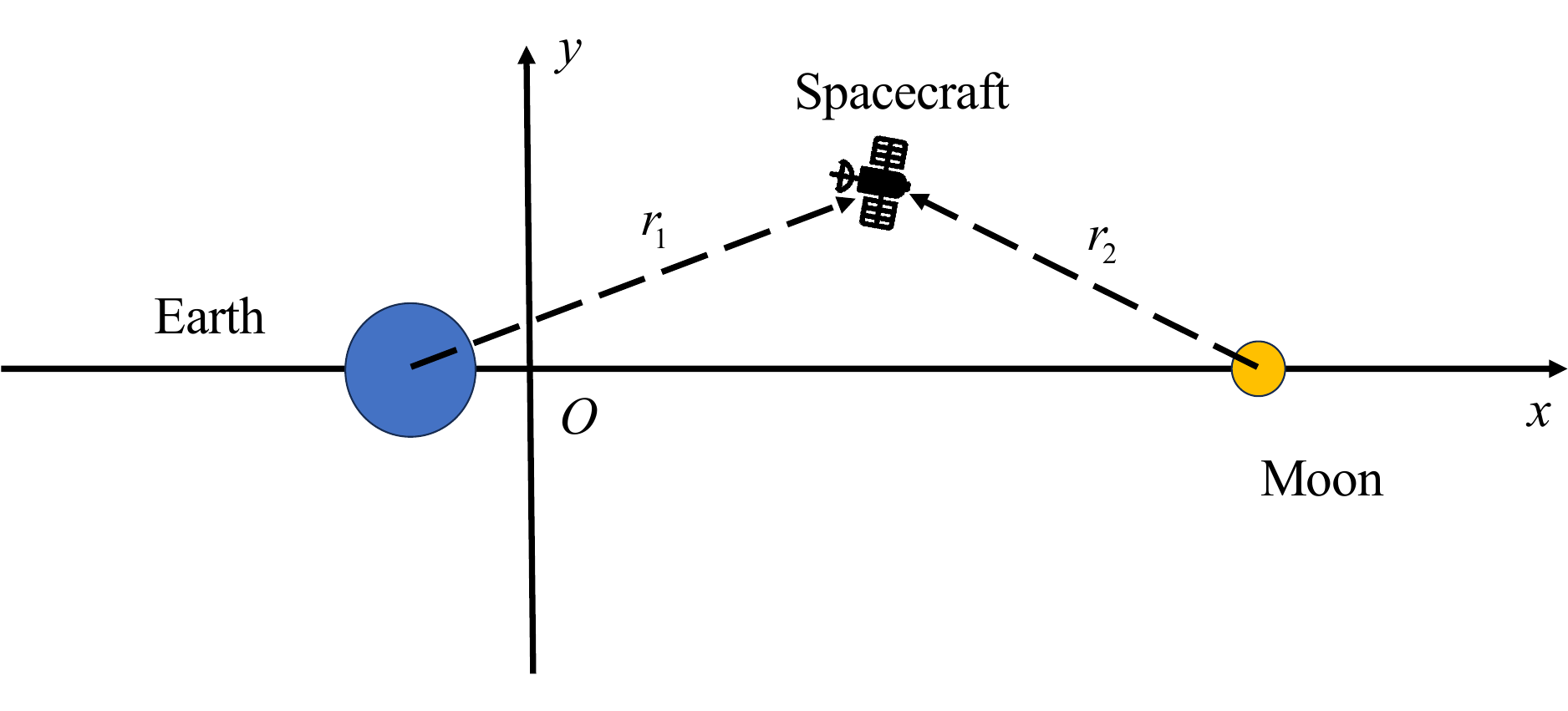}
	\caption{Earth--Moon rotating frame used in PCRTBP.}
	\label{fig:pcrtbp}
\end{figure}

The planar circular restricted three-body problem (PCRTBP) is adopted as the dynamical model in this study. In this model, the Earth and Moon move on circular orbits about their common barycenter, and the spacecraft is treated as a massless particle moving in the same plane. Figure~\ref{fig:pcrtbp} shows the barycentric rotating frame. Its $x$-axis points from Earth to Moon, and its $y$-axis completes the in-plane right-handed frame. The PCRTBP equations are normalized using the characteristic
mass, length, and time scales
$M_{\mathrm{EM}}=m_e+m_m$,
$L_{\mathrm{EM}}$, and
$t_{\mathrm{EM}}=T_{\mathrm{EM}}/(2\pi)$, respectively.
Here, $m_e$ and $m_m$ are the masses of the Earth and the Moon,
$L_{\mathrm{EM}}$ is the mean Earth-Moon distance, and $T_{\mathrm{EM}}$
is the orbital period of the Earth-Moon system. The mass parameter is $\mu=m_m/(m_e+m_m)$, giving normalized Earth and Moon masses of $1-\mu$ and $\mu$, respectively. The gravitational constant and rotating-frame angular rate are then unity. By defining the state vector as $\mathbf{X}=[x,y,v_x,v_y]^{\mathrm T}$, the equations of motion for the spacecraft can be expressed as:

\begin{equation}
    \dot{\mathbf{X}}
    =
    \begin{bmatrix}
        \dot{x}\\
        \dot{y}\\
        \dot{v}_x\\
        \dot{v}_y
    \end{bmatrix}
    =
    \begin{bmatrix}
        v_x\\
        v_y\\
        2v_y+\dfrac{\partial U}{\partial x}\\
        -2v_x+\dfrac{\partial U}{\partial y}
    \end{bmatrix},
    \label{eq:pcrtbp_eom}
\end{equation}
where $U(x,y)$ represents the pseudo-potential function of the rotating system, defined as:
\begin{equation}
\begin{aligned}
    U(x,y)&=\frac{1}{2}(x^2+y^2)+\frac{1-\mu}{r_1}+\frac{\mu}{r_2},\\
\end{aligned}
    \label{eq:pseudopotential}
\end{equation}

The terms
$r_1=\sqrt{(x+\mu)^2+y^2}$ and
$r_2=\sqrt{(x-1+\mu)^2+y^2}$
represent the normalized distances from the spacecraft to the Earth and the Moon, respectively.

Trajectory sensitivity and nonlinear targeting require the local dynamics about a nominal flow $\mathbf{X}^{*}(t)$. For the perturbed state
$\mathbf{X}(t)=\mathbf{X}^{*}(t)+\delta\mathbf{X}(t)$, the state transition matrix (STM) $\boldsymbol{\Phi}$ and second-order state transition tensor (STT) $\boldsymbol{\Psi}$ approximate the propagated deviation as
\begin{equation}
	\delta\mathbf{X}(t)
	=
	\boldsymbol{\Phi}(t,t_0)\delta\mathbf{X}(t_0)
	+
	\frac{1}{2}
	\boldsymbol{\Psi}(t,t_0)
	\delta\mathbf{X}(t_0)
	\delta\mathbf{X}(t_0),
	\label{eq:second_order_deviation_map}
\end{equation}
where $\boldsymbol{\Phi}$ is a $4\times4$ matrix and $\boldsymbol{\Psi}$ is a $4\times4\times4$ tensor.
These transition operators describe the continuous deformation of the phase space and are obtained by simultaneously integrating the following differential equations alongside the nominal state:
\begin{equation}
	\left\{
	\begin{aligned}
		\dot{\boldsymbol{\Phi}}(t,t_0)
		&=
		\mathbf{A}(t)\boldsymbol{\Phi}(t,t_0),
		\\
		\dot{\boldsymbol{\Psi}}(t,t_0)
		&=
		\mathbf{A}(t)\boldsymbol{\Psi}(t,t_0)
		+
		\mathbf{B}(t)
		\boldsymbol{\Phi}(t,t_0)
		\boldsymbol{\Phi}(t,t_0),
	\end{aligned}
	\right.
	\label{eq:stm_stt_dynamics}
\end{equation}
subject to the initial conditions
$\boldsymbol{\Phi}(t_0,t_0)=\mathbf{I}_{4\times4}$
and
$\boldsymbol{\Psi}(t_0,t_0)=\mathbf{0}_{4\times4\times4}$.
where
$\mathbf{A}(t)=\partial\mathbf{f}/\partial\mathbf{X}$
is the $4\times4$ time-varying Jacobian matrix evaluated along $\mathbf{X}^{*}(t)$, and
$\mathbf{B}(t)=\partial^{2}\mathbf{f}/\partial\mathbf{X}^{2}$
is the $4\times4\times4$ Hessian tensor of the vector field. The STM supports the sensitivity analysis, while the STT is used in the nonlinear correction targeting.

\subsection{Low-Energy Earth--Moon Transfers in PCRTBP}\label{subsec2.2}

\begin{figure}[htbp]
\centering
    \includegraphics[width=0.5\textwidth]{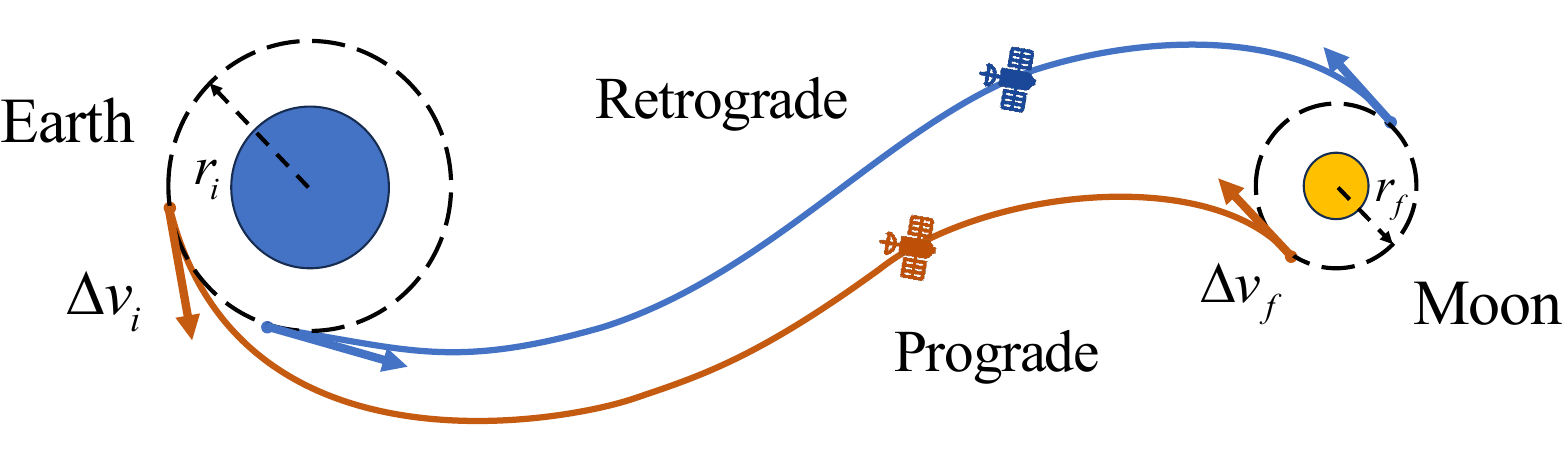}
    \caption{Departure and arrival impulses for two-impulse Earth--Moon transfer.}
    \label{fig:transfer_geometry}
\end{figure}

In this study, we particularly focus on the two-impulse transfers from a circular low Earth orbit to a circular low lunar orbit in PCRTBP. The departure and arrival radii are $r_i$ and $r_f$. A tangential impulse $\Delta v_i$ starts the transfer, and $\Delta v_f$ inserts the spacecraft into the lunar orbit (Fig.~\ref{fig:transfer_geometry}). Prograde and retrograde cases are distinguished by the direction of lunar arrival. 

The boundary constraints are
\begin{equation}
    \boldsymbol{\psi}_i(\mathbf{X}_i)
    =
    \begin{bmatrix}
        (x_i+\mu)^2+y_i^2-r_i^2\\
        (x_i+\mu)(\dot{x}_i-y_i)+y_i(\dot{y}_i+x_i+\mu)
    \end{bmatrix}
    =\mathbf{0}.
    \label{eq:departure_constraint}
\end{equation}
and
\begin{equation}
    \boldsymbol{\psi}_f(\mathbf{X}_f)
    =
    \begin{bmatrix}
        (x_f-1+\mu)^2+y_f^2-r_f^2\\
        (x_f-1+\mu)(\dot{x}_f-y_f)+y_f(\dot{y}_f+x_f-1+\mu)
    \end{bmatrix}
    =\mathbf{0}.
    \label{eq:arrival_constraint}
\end{equation}

The terminal state is obtained by integrating Eq.~\eqref{eq:pcrtbp_eom}. A constrained nonlinear program then minimizes
\begin{equation}
    \Delta v=\Delta v_i+\Delta v_f = \sqrt{(\dot{x}_i-y_i)^2+(\dot{y}_i+x_i+\mu)^2}
    -\sqrt{\frac{1-\mu}{r_i}}
    +\sqrt{(\dot{x}_f-y_f)^2+(\dot{y}_f+x_f+\mu-1)^2}
    -\sqrt{\frac{\mu}{r_f}} .
    \label{eq:nominal_cost}
\end{equation}

\begin{figure}[htbp]
\centering
    \includegraphics[width=0.8\textwidth]{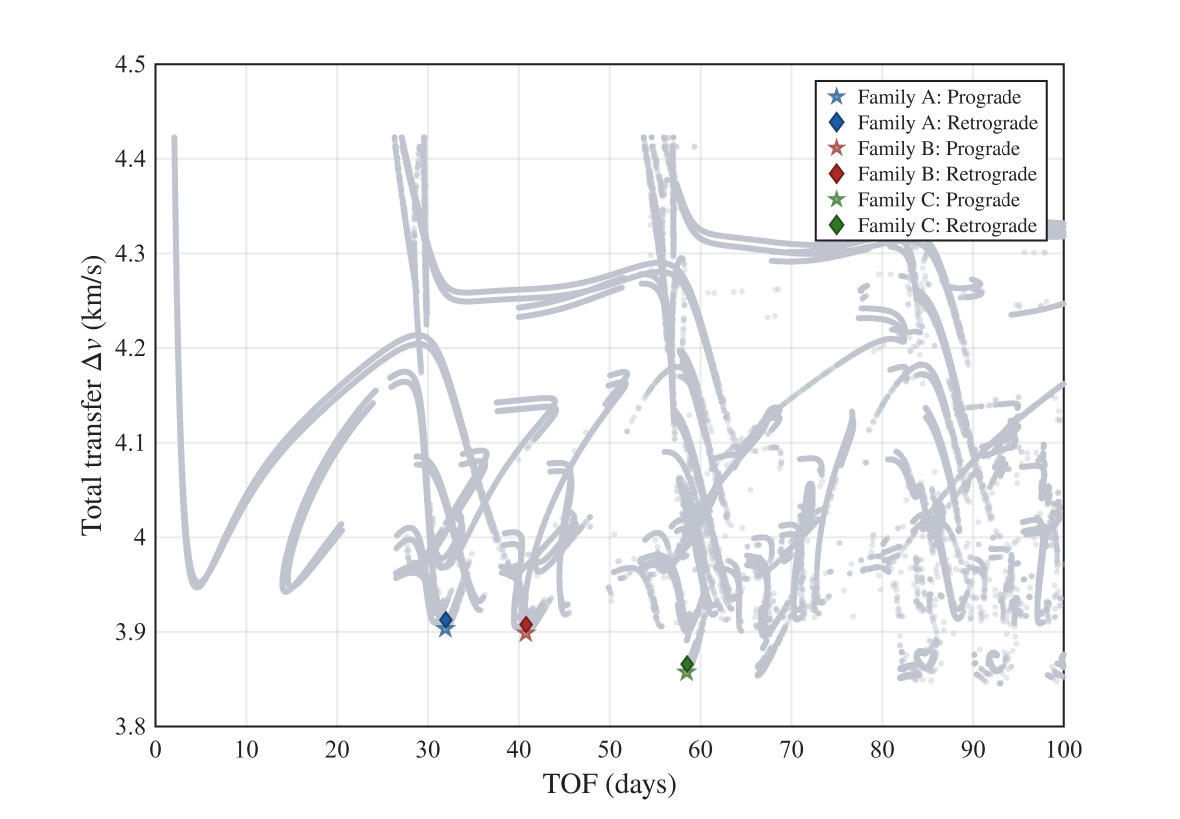}
    \caption{Transfer time and total impulse of the computed solution set. Selected cases T1--T6 are highlighted.}
    \label{fig:solution_set}
\end{figure}

\begin{figure}[htbp]
    \centering
    \includegraphics[width=0.9\textwidth]{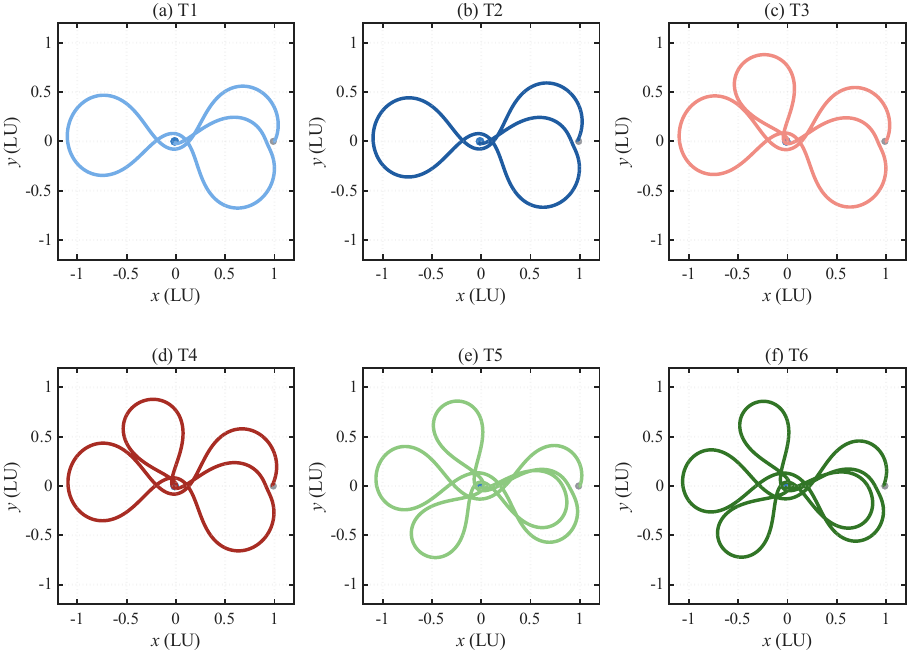}
    \caption{Reference transfers T1--T6 in the Earth--Moon rotating frame.}
    \label{fig:six_trajectories}
\end{figure}

The design variables are the departure state, arrival state, and time of flight. Figure~\ref{fig:solution_set} shows the solution set from a 167-km-altitude Earth parking orbit to a 100-km-altitude lunar orbit. Six typical local optima, T1--T6, are selected to span distinct encounter sequences and both lunar-arrival directions. Their flight times range from 30 to 60 days, and their total impulse is below \SI{3.92}{\kilo\metre\per\second}. Their geometry and parameters are given in Fig.~\ref{fig:six_trajectories} and Table~\ref{tab:nominal_transfers}.

\begin{table}[!hbt]
	\centering
	\caption{Nominal parameters of the reference transfers.}
	\label{tab:nominal_transfers}
	\begin{tabular}{clcccc}
		\toprule
		Case & Class  & TOF (day) & $\Delta v_i$ (km/s) & $\Delta v_f$ (km/s) & Total $\Delta v$ (km/s) \\
		\midrule
		$\mathrm{T1}$ & Retrograde   & 31.92 & 3.1284 & 0.7751 & 3.9035 \\
		$\mathrm{T2}$ & Prograde & 31.95 & 3.1283 & 0.7843 & 3.9126 \\
		$\mathrm{T3}$ & Retrograde   & 40.77 & 3.1275 & 0.7711 & 3.8987 \\
		$\mathrm{T4}$ & Prograde & 40.79 & 3.1274 & 0.7803 & 3.9078 \\
		$\mathrm{T5}$ & Retrograde   & 58.45 & 3.1201 & 0.7373 & 3.8574 \\
		$\mathrm{T6}$ & Prograde & 58.53 & 3.1199 & 0.7461 & 3.8661 \\
		\bottomrule
	\end{tabular}
\end{table}

\subsection{Stochastic Mid-Course Correction Formulation}\label{subsec2.3}
\begin{figure}[htbp]
\centering
    \includegraphics[width=0.8\textwidth]{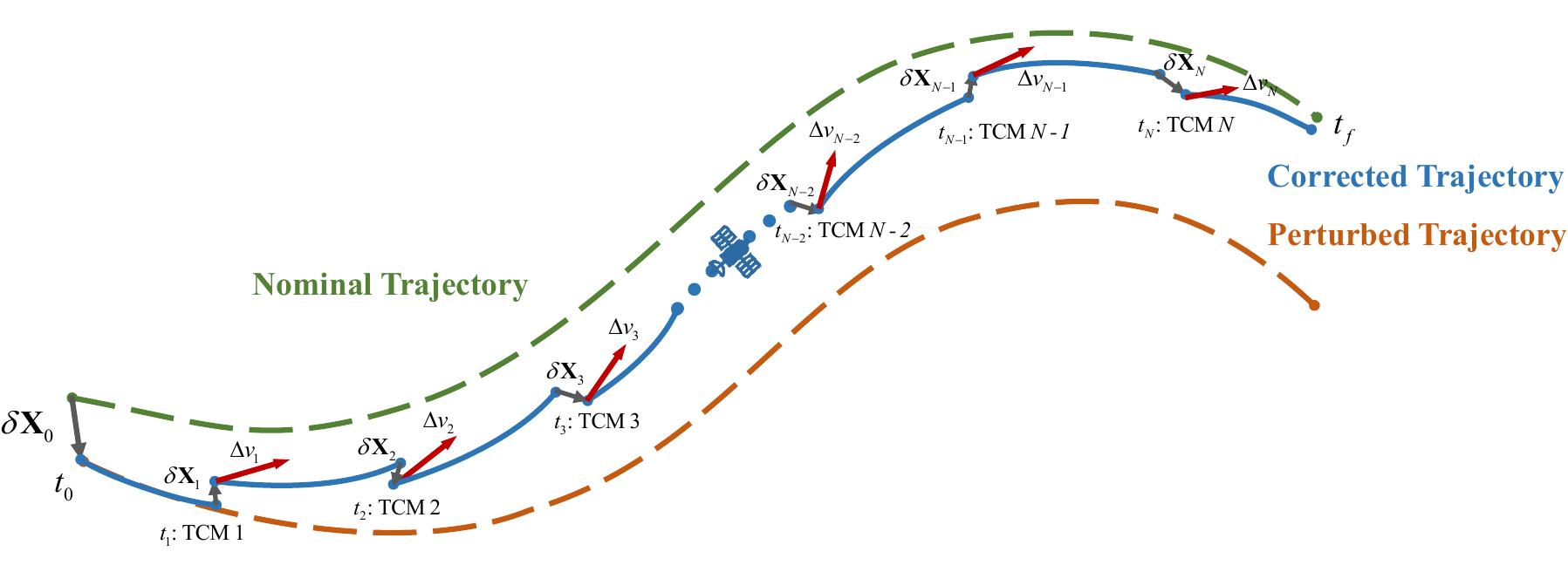}
    \caption{Stochastic correction model with navigation and maneuver-execution errors.}
    \label{fig:correction_model}
\end{figure}

A correction strategy must determine the maneuver count $N_{\mathrm{TCM}}$ and the correction epochs $\mathcal{T}=\{t_k\}_{k=1}^{N_{\mathrm{TCM}}}$. Let $\boldsymbol{\xi}$ collect the injection, navigation, and maneuver-execution errors. Corrections reduce the accumulated state deviation, while each executed command introduces a new error. Figure~\ref{fig:correction_model} shows this estimation--correction--propagation cycle.

For a given architecture $(N_{\mathrm{TCM}},\mathcal{T})$, the total correction cost is defined as
\begin{equation}
    D\left(N_{\mathrm{TCM}},\mathcal{T};\boldsymbol{\xi}\right)
    =\sum_{k=1}^{N_{\mathrm{TCM}}}\left\|\Delta\mathbf{v}_k\right\|.
    \label{eq:correction_cost}
\end{equation}

The correction strategy is therefore formulated as a stochastic optimization problem. The objective is to reduce the correction demand while maintaining a prescribed probability of successful trajectory recovery:
\begin{equation}
\begin{aligned}
	\min_{N_{\mathrm{TCM}},\,\mathcal{T}}
	\quad &
	\mathcal{R}_D
	=
	\mathbb{E}[D\mid\mathcal{E}]
	+
	\kappa
	\sqrt{\operatorname{Var}(D\mid\mathcal{E})},
	\\
	\text{s.t.}\quad &
	\mathbb{P}
	\left[
	\mathcal{E}
	(N_{\mathrm{TCM}},\mathcal{T};\boldsymbol{\xi})
	\right]
	\geq p_{\mathrm{req}},
	\\
	&
	t_0<t_1<t_2<\cdots<t_{N_{\mathrm{TCM}}}<t_f .
\end{aligned}
\label{eq:general_stochastic_problem}
\end{equation}
where, $\mathcal{E}$ denotes the joint mission-success event, $p_{\mathrm{req}}$ is the required reliability, and $\kappa$ controls the penalty on correction-cost dispersion. The criteria for joint mission success are described in detail in Section~\ref{sec5}. The conditional fuel metric is evaluated only over successful realizations, while unsuccessful realizations are accounted for through the reliability constraint. In the proposed method, the discrete maneuver architecture is determined from the local sensitivity structure in Section~\ref{sec3}, while the continuous correction epochs are optimized in Section~\ref{sec4}.

\section{Local Error Amplification and Sensitivity-Guided Segmentation}\label{sec3}

This section derives the instantaneous local Lyapunov exponent (LLE) from the variational dynamics and uses it to identify sensitive passages and segment the nominal trajectory. Numerical simulations are conducted to verify the relationship between LLE maxima and finite-error amplification.

\subsection{Instantaneous Local Lyapunov Exponent}\label{subsec3.1}
Let $\mathbf{X}^{*}(t)$ be a nominal PCRTBP trajectory and let $\delta\mathbf{X}(t)$ be an infinitesimal state deviation. Linearization of the equations of motion gives
\begin{equation}
    \delta\dot{\mathbf{X}}(t)
    =\mathbf{A}(t)\delta\mathbf{X}(t),
    \qquad
    \mathbf{A}(t)
    =\left.
    \frac{\partial\mathbf{f}}{\partial\mathbf{X}}
    \right|_{\mathbf{X}^{*}(t)}.
    \label{eq:linear_variational_dynamics}
\end{equation}

Correction scheduling requires a local measure of error amplification along the nominal trajectory. Therefore, the LLE is obtained from the zero-horizon limit of the forward finite-time Lyapunov exponent (FTLE)~\cite{nolan2020ftle}. For \(\Delta t>0\), let \(\boldsymbol{\Phi}(t+\Delta t,t)\) denote the STM associated with Eq.~\eqref{eq:linear_variational_dynamics}. In the adopted normalized state coordinates and Euclidean norm, the right Cauchy--Green tensor and the FTLE are
\begin{equation}
    \mathbf{C}(t,\Delta t)
    =\boldsymbol{\Phi}^{\mathrm T}(t+\Delta t,t)
    \boldsymbol{\Phi}(t+\Delta t,t),
    \qquad
    \sigma(t,\Delta t)
    =\frac{1}{2\Delta t}
    \ln\lambda_{\max}\!\left[\mathbf{C}(t,\Delta t)\right].
    \label{eq:ftle_definition}
\end{equation}

Since \(\dot{\boldsymbol{\Phi}}=\mathbf{A}\boldsymbol{\Phi}\) and \(\boldsymbol{\Phi}(t,t)=\mathbf{I}\), a short-time expansion gives
\begin{equation}
\begin{aligned}
    \boldsymbol{\Phi}(t+\Delta t,t)
    &=\mathbf{I}+\mathbf{A}(t)\Delta t+\mathcal{O}(\Delta t^2),\\
    \mathbf{C}(t,\Delta t)
    &=\mathbf{I}+2\mathbf{S}(t)\Delta t+\mathcal{O}(\Delta t^2),
\end{aligned}
    \label{eq:short_time_cg}
\end{equation}
where \(\mathbf{S}(t)=\tfrac{1}{2}[\mathbf{A}(t)+\mathbf{A}^{\mathrm T}(t)]\). Because \(\mathbf{S}\) is symmetric and \(\Delta t>0\), its largest eigenvalue determines the leading stretching term:
\begin{equation}
    \lambda_{\max}\!\left[\mathbf{C}(t,\Delta t)\right]
    =1+2\lambda_{\max}\!\left[\mathbf{S}(t)\right]\Delta t
    +\mathcal{O}(\Delta t^2).
    \label{eq:short_time_cg_eigenvalue}
\end{equation}
Substituting this expression into Eq.~\eqref{eq:ftle_definition} and using \(\ln(1+z)=z+\mathcal{O}(z^2)\) yields
\begin{equation}
    \Lambda_{\mathrm{LLE}}(t)
    =\lim_{\Delta t\rightarrow0^+}\sigma(t,\Delta t)
    =\lambda_{\max}\!\left[\mathbf{S}(t)\right].
    \label{eq:lle_general}
\end{equation}
Thus, the LLE measures the maximum instantaneous logarithmic stretching rate under the chosen state normalization. The corresponding eigenvector of \(\mathbf{S}\) identifies the perturbation direction that realizes this rate.

The PCRTBP structure permits a closed-form expression for $\Lambda_{\mathrm{LLE}}$. Define
\begin{equation}
    \mathbf{H}_U
    =\begin{bmatrix}
        U_{xx} & U_{xy}\\
        U_{xy} & U_{yy}
    \end{bmatrix},
    \qquad
    \mathbf{R}
    =\begin{bmatrix}
        0 & 2\\
        -2 & 0
    \end{bmatrix}.
    \label{eq:pcrtbp_blocks}
\end{equation}

Using Eq.~\eqref{eq:pcrtbp_eom}, the variational matrix and its symmetric part become
\begin{equation}
    \mathbf{A}
    =\begin{bmatrix}
        \mathbf{0} & \mathbf{I}\\
        \mathbf{H}_U & \mathbf{R}
    \end{bmatrix},
    \qquad
    \mathbf{S}
    =\begin{bmatrix}
        \mathbf{0} & \mathbf{K}\\
        \mathbf{K} & \mathbf{0}
    \end{bmatrix},
    \qquad
    \mathbf{K}=\frac{1}{2}(\mathbf{I}+\mathbf{H}_U).
    \label{eq:symmetric_block_form}
\end{equation}

The Coriolis matrix $\mathbf{R}$ is skew-symmetric and does not contribute to $\mathbf{S}$. The two eigenvalues of $\mathbf{K}$ are
\begin{equation}
    \kappa_{\pm}
    =\frac{1}{2}
    +\frac{U_{xx}+U_{yy}
    \pm\sqrt{(U_{xx}-U_{yy})^2+4U_{xy}^2}}{4}.
    \label{eq:kappa_eigenvalues}
\end{equation}

Since $\mathbf{K}$ is symmetric, the eigenvalues of $\mathbf{S}$ are $\{\pm\kappa_{+},\pm\kappa_{-}\}$. In the PCRTBP,
\begin{equation}
    U_{xx}+U_{yy}
    =2+\frac{1-\mu}{r_1^3}+\frac{\mu}{r_2^3}>0.
    \label{eq:pseudopotential_trace}
\end{equation}

Therefore, $\kappa_{+}+\kappa_{-}>0$, which implies $\kappa_{+}>-\kappa_{-}$. The largest eigenvalue of $\mathbf{S}$ is consequently $\kappa_{+}$, and
\begin{equation}
    \Lambda_{\mathrm{LLE}}(x,y)
    =\frac{1}{2}
    +\frac{U_{xx}+U_{yy}
    +\sqrt{(U_{xx}-U_{yy})^2+4U_{xy}^2}}{4}.
    \label{eq:lle_closed}
\end{equation}

Equation~\eqref{eq:lle_closed} depends only on the current position and the local curvature of the pseudo-potential. Its connection with near-primary passages can be made explicit. For primary $j$, define $\boldsymbol{\rho}_j=r_j\mathbf{e}_{r,j}$ and let $m_j$ denote its normalized mass. The pseudo-potential Hessian is
\begin{equation}
    \mathbf{H}_U
    =\mathbf{I}
    +\sum_{j\in\{\mathrm E,\mathrm M\}}
    \frac{m_j}{r_j^3}
    \left(3\mathbf{e}_{r,j}\mathbf{e}_{r,j}^{\mathrm T}
    -\mathbf{I}\right),
    \qquad
    m_{\mathrm E}=1-\mu,
    \quad
    m_{\mathrm M}=\mu.
    \label{eq:pseudopotential_hessian_primary_form}
\end{equation}

Near primary $j$, let $\gamma_j=m_j/r_j^3$ denote the dominant gravity-gradient term. In the local radial--tangential basis,
\begin{equation}
\begin{aligned}
    \lambda_r(\mathbf{H}_U)
    &=1+2\gamma_j+\mathcal{O}(1),\\
    \lambda_t(\mathbf{H}_U)
    &=1-\gamma_j+\mathcal{O}(1).
\end{aligned}
    \label{eq:near_primary_hessian_eigenvalues}
\end{equation}

The corresponding eigenvalues of $\mathbf{K}$ are $1+\gamma_j+\mathcal{O}(1)$ and $1-\gamma_j/2+\mathcal{O}(1)$. Hence,
\begin{equation}
    \Lambda_{\mathrm{LLE}}
    =1+\frac{m_j}{r_j^3}+\mathcal{O}(1),
    \quad
    \frac{\partial\Lambda_{\mathrm{LLE}}}{\partial r_j}
    \simeq-\frac{3m_j}{r_j^4}<0.
    \label{eq:near_primary_lle_scaling}
\end{equation}

The LLE therefore rises approximately with the inverse cube of the distance to the nearby primary. The radial curvature grows as $2\gamma_j$, while the tangential curvature changes as $-\gamma_j$. An initially compact uncertainty set should consequently become elongated during a periapsis passage. The LLE remains a local upper growth rate rather than a finite-time covariance prediction. Equations~\eqref{eq:lle_general} and \eqref{eq:near_primary_lle_scaling} identify where strong stretching is expected, but they do not determine how much a finite error grows over a prescribed interval. Both points are examined numerically below.

\subsection{Local Sensitivity and Finite-Time Error Amplification}\label{subsec3.2}

\begin{figure}[!htbp]
    \centering
    \captionsetup[subfigure]{font=footnotesize,skip=1pt}
    \begin{subfigure}[b]{0.31\textwidth}
        \centering
        \includegraphics[width=\textwidth]{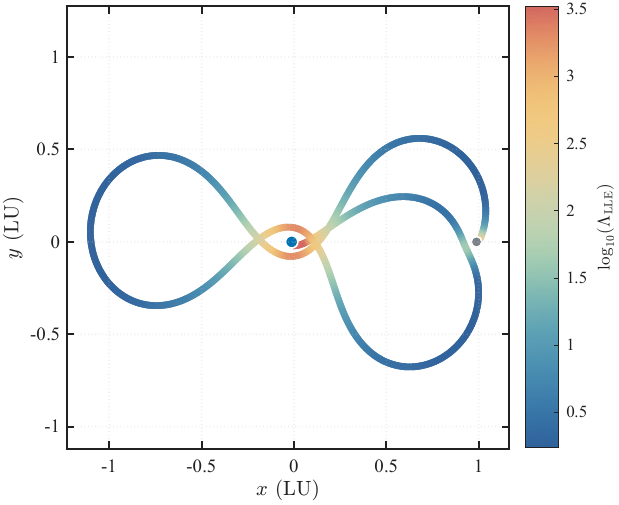}
        \caption{T1}
    \end{subfigure}
    \hfill
    \begin{subfigure}[b]{0.31\textwidth}
        \centering
        \includegraphics[width=\textwidth]{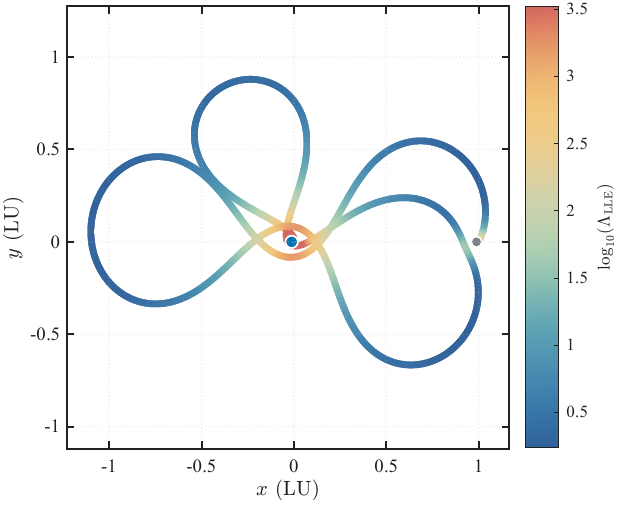}
        \caption{T3}
    \end{subfigure}
    \hfill
    \begin{subfigure}[b]{0.31\textwidth}
        \centering
        \includegraphics[width=\textwidth]{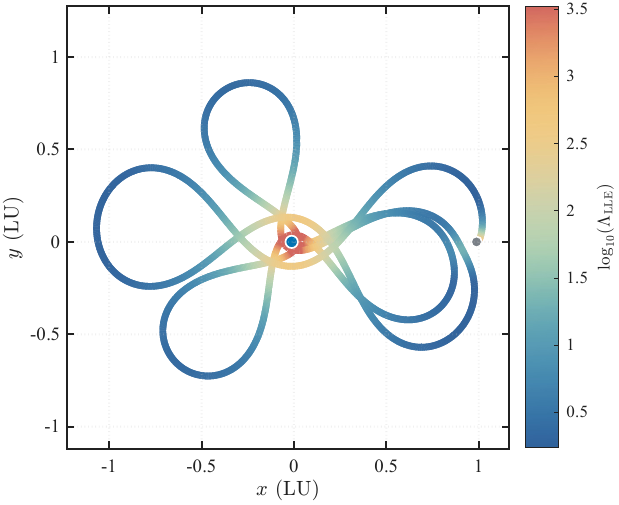}
        \caption{T5}
    \end{subfigure}
    \\[0.15em]
    \begin{subfigure}[b]{0.31\textwidth}
        \centering
        \includegraphics[width=\textwidth]{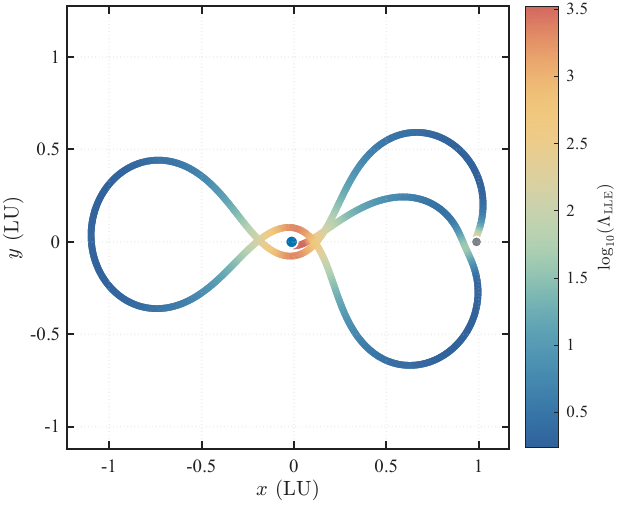}
        \caption{T2}
    \end{subfigure}
    \hfill
    \begin{subfigure}[b]{0.31\textwidth}
        \centering
        \includegraphics[width=\textwidth]{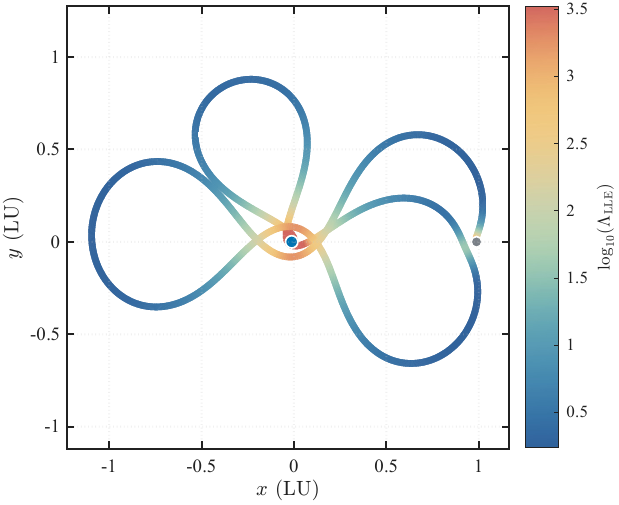}
        \caption{T4}
    \end{subfigure}
    \hfill
    \begin{subfigure}[b]{0.31\textwidth}
        \centering
        \includegraphics[width=\textwidth]{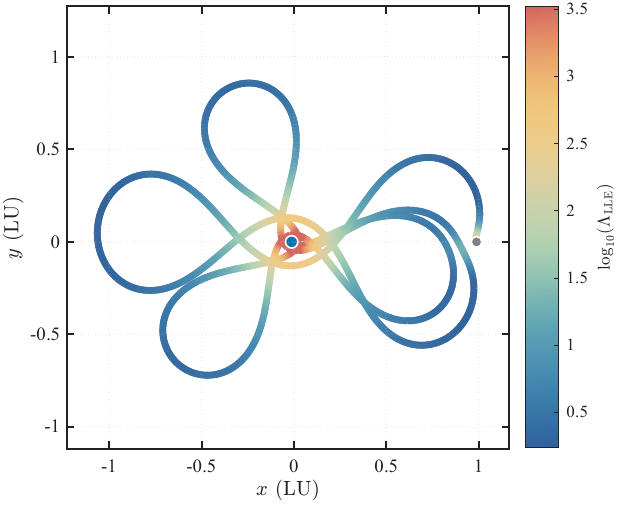}
        \caption{T6}
    \end{subfigure}

    \caption{Spatial distribution of the instantaneous LLE along T1--T6. }
    \label{fig.LLE_spatial}
\end{figure}

\begin{figure}[!htbp]
    \centering
    \captionsetup[subfigure]{font=footnotesize,skip=1pt}
    \begin{subfigure}[b]{0.31\textwidth}
        \centering
        \includegraphics[width=\textwidth]{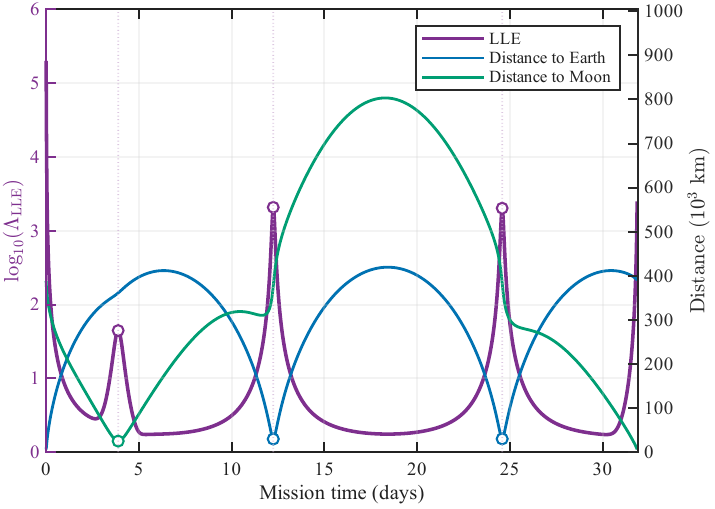}
        \caption{T1}
    \end{subfigure}
    \hfill
    \begin{subfigure}[b]{0.31\textwidth}
        \centering
        \includegraphics[width=\textwidth]{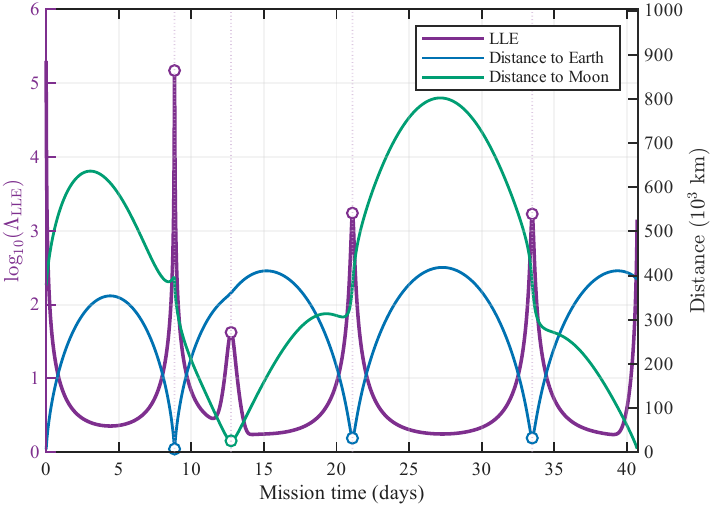}
        \caption{T3}
    \end{subfigure}
    \hfill
    \begin{subfigure}[b]{0.31\textwidth}
        \centering
        \includegraphics[width=\textwidth]{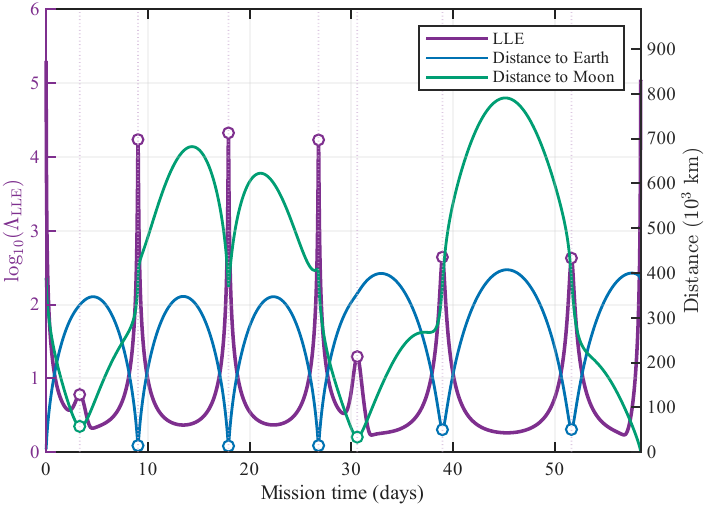}
        \caption{T5}
    \end{subfigure}
    \\[0.15em]
    \begin{subfigure}[b]{0.31\textwidth}
        \centering
        \includegraphics[width=\textwidth]{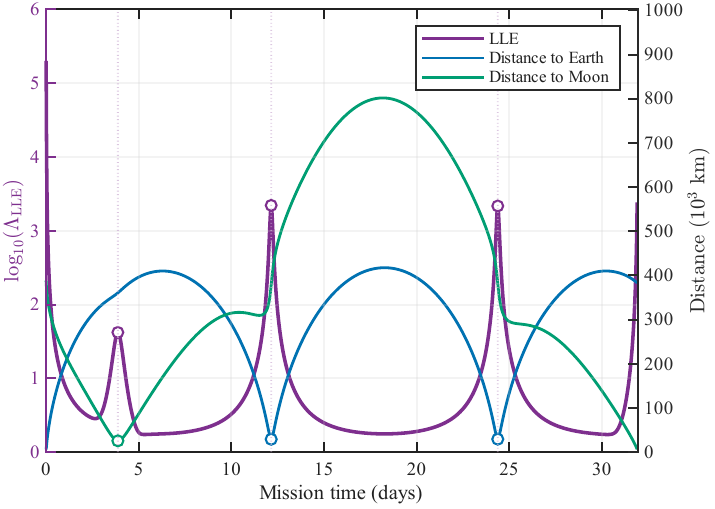}
        \caption{T2}
    \end{subfigure}
    \hfill
    \begin{subfigure}[b]{0.31\textwidth}
        \centering
        \includegraphics[width=\textwidth]{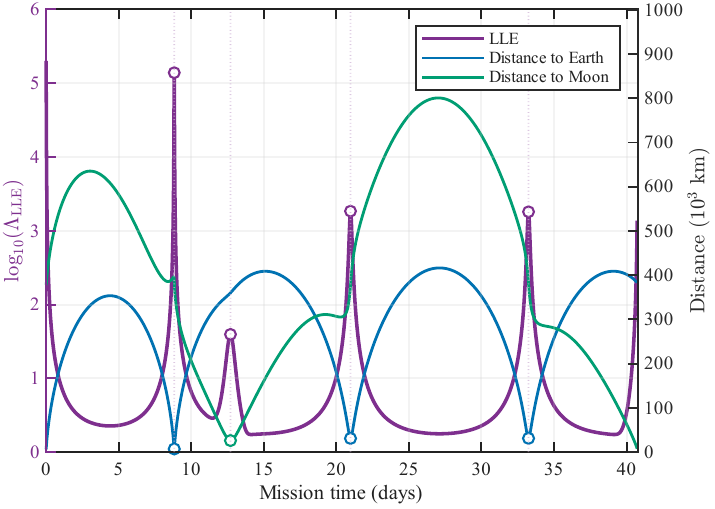}
        \caption{T4}
    \end{subfigure}
    \hfill
    \begin{subfigure}[b]{0.31\textwidth}
        \centering
        \includegraphics[width=\textwidth]{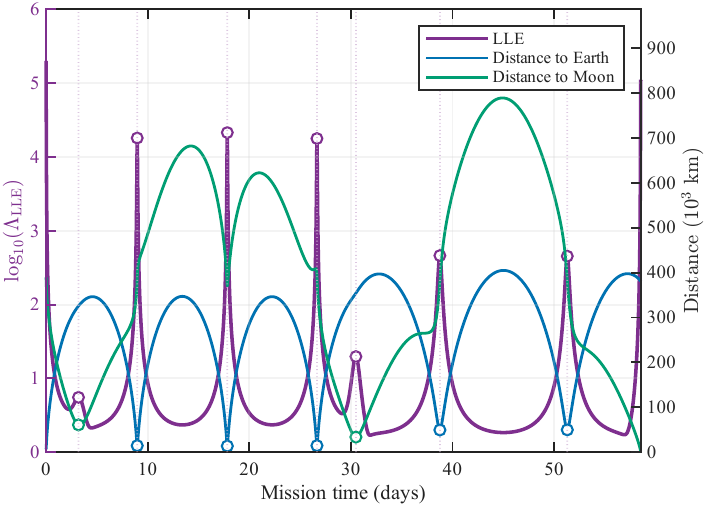}
        \caption{T6}
    \end{subfigure}

    \caption{LLE and primary-distance for T1--T6. Purple circles mark the  LLE local maxima. The corresponding blue or green circles mark the distance to the nearer primary at the same epoch.}
    \label{fig.LLE_seg}
\end{figure}

Figure.~\ref{fig.LLE_spatial} and \ref{fig.LLE_seg} show that the LLE increases sharply near Earth perigees and lunar perilunes, while remaining much lower over most of the intervening arcs. For all six transfers, the strongest interior LLE peak occurs at an Earth perigee (Table~\ref{tab:sensitivity_structure}). These results agree with the near-primary relationship derived in Eq.~\eqref{eq:near_primary_lle_scaling}. As the spacecraft approaches a primary, the local gravity gradient grows approximately as $m_j/r_j^3$, which increases the largest eigenvalue of the symmetric variational dynamics and creates a narrow high-LLE region. The concentration of high-LLE regions around close approaches therefore supports the analytical connection between local gravity gradients and instantaneous dynamical sensitivity. The recurring Earth and lunar peaks identify distinct sensitive passages along each transfer, providing a dynamical basis for the subsequent trajectory segmentation.

\begin{figure}[htbp]
    \centering
    \includegraphics[width=0.55\textwidth]{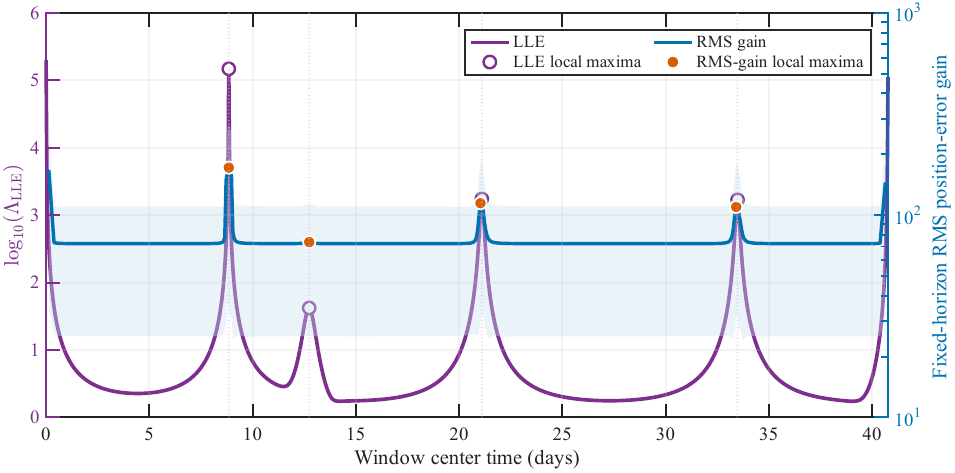}
    \caption{Finite-time comparison of the instantaneous LLE and RMS position-error gain on T3.}
    \label{fig:lle_local_mc_validation}
\end{figure}

\begin{table}[htbp]
\centering
\caption{Sensitivity-event structure and the strongest event along each nominal transfer. Perigee radius is measured from the Earth's center.}
\label{tab:sensitivity_structure}
\small
\begin{tabular}{cccccc}
\toprule
Case & Interior peaks & Segments & Strongest $\log_{10}\Lambda_{\mathrm{LLE}}$ & Location & Perigee radius (km)\\
\midrule
T1 & 3 & 4 & 3.32 & Earth perigee & 29982\\
T2 & 3 & 4 & 3.34 & Earth perigee & 29386\\
T3 & 4 & 5 & 5.17 & Earth perigee & 7233\\
T4 & 4 & 5 & 5.14 & Earth perigee & 7396\\
T5 & 7 & 8 & 4.33 & Earth perigee & 13827\\
T6 & 7 & 8 & 4.33 & Earth perigee & 13786\\
\bottomrule
\end{tabular}
\end{table}
\FloatBarrier
Finite-time error amplification is examined on T3 by restarting the same perturbation ensemble at successive epochs and propagating each ensemble over a fixed horizon. The resulting position-error amplification is measured by
\begin{equation}
    G_{\mathrm{RMS}}(t_j)=
    \frac{
    	\sqrt{\dfrac{1}{N_{\mathrm{MC}}}
    		\sum_{q=1}^{N_{\mathrm{MC}}}
    		\left|\delta\mathbf r^{(q)}(t_j+T_h/2)\right|^2}
    }{
    	\sqrt{\dfrac{1}{N_{\mathrm{MC}}}
    		\sum_{q=1}^{N_{\mathrm{MC}}}
    		\left|\delta\mathbf r^{(q)}(t_j-T_h/2)\right|^2}
    }.
\end{equation}
where $N_{\mathrm{MC}}$ is the number of perturbation samples and $T_h$ is the fixed propagation horizon. The ratio compares the RMS position error after propagation with that at the restart epoch, so $G_{\mathrm{RMS}}>1$ indicates finite-time error amplification. Because the same perturbation ensemble and horizon are used at every restart, variations in $G_{\mathrm{RMS}}$ primarily reflect changes in the local dynamics. Each interior LLE maximum in Fig.~\ref{fig:lle_local_mc_validation} is accompanied by a rise in finite-time error gain. The gain is broader and can be slightly displaced because it accumulates the dynamics over a finite interval, while the LLE is instantaneous.

\begin{figure}[htbp]
    \centering
    \includegraphics[width=0.7\textwidth]{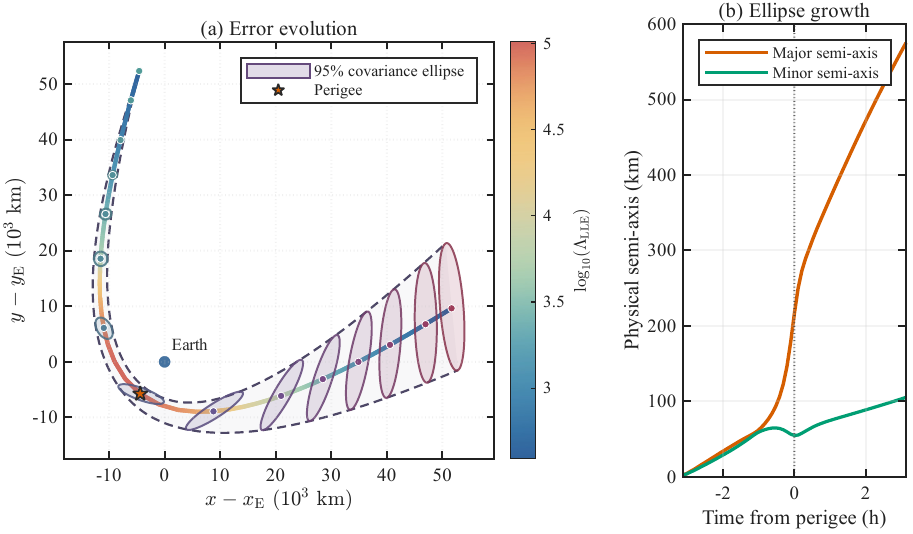}
    \caption{T3 error evolution through the strongest Earth perigee. (a) Local trajectory, 95\% position-covariance ellipses, and their envelope; (b) 95\% covariance semi-axis growth.}
    \label{fig:perigee_uncertainty_ellipses}
\end{figure}

Figure~\ref{fig:perigee_uncertainty_ellipses} further illustrates the directional nature of the local error amplification. In Fig.~\ref{fig:perigee_uncertainty_ellipses} (a), the initially compact error cloud becomes strongly elongated through the highest-LLE Earth perigee, predominantly along the local radial direction, while the transverse deformation remains much weaker. This behavior follows the near-primary scaling in Eq.~\eqref{eq:near_primary_hessian_eigenvalues}. As $r_E$ decreases, $\gamma_E=m_E/r_E^3$ rises rapidly, increasing the radial curvature $1+2\gamma_E$ much more strongly than the tangential curvature $1-\gamma_E$. Consistently, Fig.~\ref{fig:perigee_uncertainty_ellipses} (b) shows a sharp increase in the covariance major semi-axis around perigee, while the minor semi-axis varies only moderately. These results show that the LLE identifies where finite errors begin to grow rapidly.

\FloatBarrier

\subsection{Sensitivity-Guided Segmentation}\label{subsec3.3}
Based on the preceding analysis, the trajectory is segmented at the internal local maxima of the LLE. Let $\{\tau_i\}_{i=1}^{N_{\mathrm p}}$ denote the ordered interior maxima of the instantaneous LLE history. Departure and arrival complete the boundary set
\begin{equation}
    \mathcal{B}
    =\left\{\tau_0,\tau_1,\ldots,
    \tau_{N_{\mathrm p}},\tau_{N_{\mathrm p}+1}\right\},
    \qquad
    \tau_0=t_0,
    \qquad
    \tau_{N_{\mathrm p}+1}=t_f.
    \label{eq:segment_boundaries}
\end{equation}
The transfer is divided into
\begin{equation}
    \mathcal{S}_i=[\tau_{i-1},\tau_i],
    \qquad
    i=1,\ldots,N_{\mathrm{seg}},
    \qquad
    N_{\mathrm{seg}}=N_{\mathrm p}+1.
    \label{eq:lle_segments}
\end{equation}
An LLE maxima marks the centre of a short sensitive passage, typically associated with an Earth perigee or lunar perilune. Maintaining sufficient state accuracy at these sensitivity peaks is therefore critical to correction efficiency. If a large deviation enters the high-LLE region, the strong local stretching shown in Fig.~\ref{fig:perigee_uncertainty_ellipses} can rapidly amplify the error and increase the downstream recovery demand. An upstream correction is thus used to reduce the deviation before the sensitive passage, while a subsequent navigation update and correction remove the stretched residual before it is carried toward the next sensitive region. The exact correction epochs are determined by the stochastic design in Section~\ref{sec4}.

\begin{figure}[htbp]
    \centering
    \includegraphics[width=0.70\textwidth]{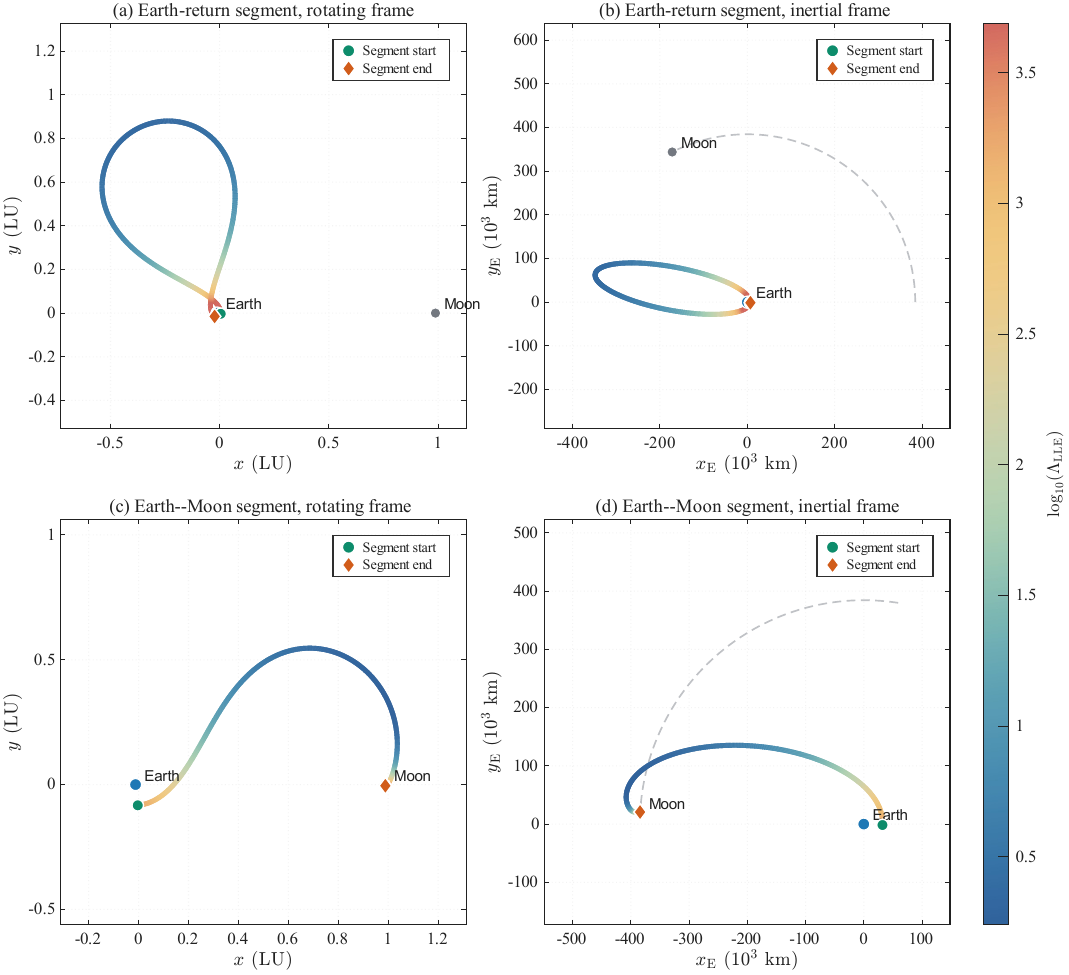}
    \caption{Representative LLE-based segments on T3, colored by the instantaneous LLE. The dashed curves in the inertial-frame panels show lunar motion during each segment.}
    \label{fig:representative_segments}
\end{figure}

The two basic segment geometries generated by the LLE-based segmentation is shown in Fig.~\ref{fig:representative_segments}. After segmentation, each low-energy transfer can be viewed as a sequence of these two types of local arcs: an Earth-centered perigee-to-perigee return arc and a direct Earth--Moon or Moon--Earth transfer arc. The segment boundaries are located at the LLE maxima, while the interior of each segment generally remains less sensitive. This decomposition converts a complex multi-encounter transfer into a sequence of simpler local targeting problems. For each segment, the correction design only needs to maintain sufficient accuracy at its downstream boundary, which simplifies the targeting and guidance formulation. Moreover, the relatively low sensitivity over most of the segment interior allows conventional mid-course correction concepts to be used locally for correction-epoch design. This segmentation provides the dynamical basis for the maneuver architecture and timing design in Section~\ref{sec4}. Notably, local Lyapunov information has previously been incorporated into differential-correction procedures for cislunar trajectory targeting~\cite{blumenthal2023lle}. In that application, the local exponent is used within the targeting process. The present study uses the same type of local sensitivity information for a different purpose: to identify sensitivity events along the complete low-energy transfer and to construct the subsequent segmentation and correction architecture.

\FloatBarrier

\section{LLE-Guided Segmented Mid-Course Correction Design and Optimization}\label{sec4}
This section constructs the correction architecture from the LLE-based segmentation and introduces a reduced timing parameterization. A second-order STT--LM targeting procedure is presented to compute the correction impulses. The correction schedule is then optimized through a chance-constrained formulation and subjected to independent reliability validation.

\subsection{Maneuver Architecture from LLE Segmentation}\label{subsec4.1}
\begin{figure}[htbp]
    \centering
    \includegraphics[width=0.6\textwidth]{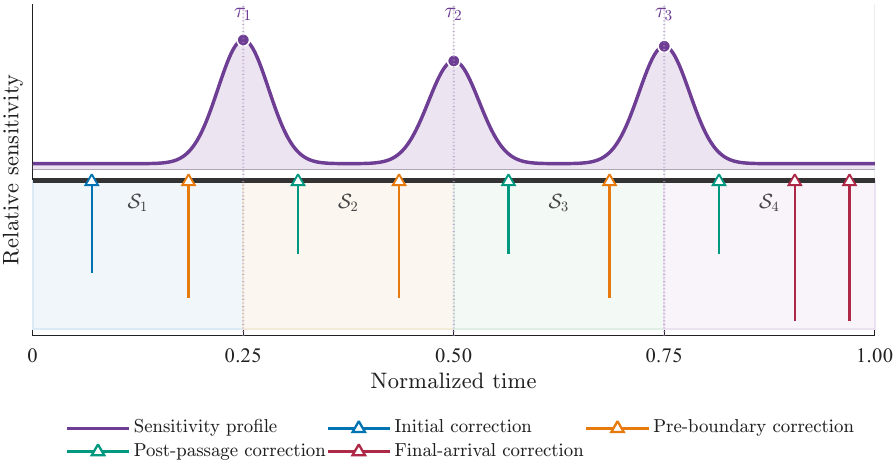}
    \caption{Schematic correction placement relative to the LLE-guided segments.}
    \label{fig:maneuver_architecture}
\end{figure}

The stochastic problem in Eq.~\eqref{eq:general_stochastic_problem} contains a discrete decision, $N_{\mathrm{TCM}}$, and a continuous schedule, $\mathcal T$. Directly searching over both would require repeated nonlinear targeting and stochastic evaluation for a large number of candidate maneuver architectures. The segmentation method introduced above allows us to first determine the number of correction maneuvers and then optimize their timing. Every correction in segment $\mathcal S_i$ targets the nominal position at its downstream boundary $\tau_i$. The correction sequence is constructed around the LLE-guided segment boundaries, as shown in Fig~\ref{fig:maneuver_architecture}. The first correction is performed shortly after injection to remove the initial delivery error. For each interior boundary $\tau_i$, corresponding to a local LLE maximum, one correction is placed upstream to ensure sufficient accuracy before entering the sensitive region, and another is placed downstream to remove the residual amplified through the passage. Finally, two corrections are assigned before lunar orbit insertion to refine the arrival state. The total number of correction maneuvers is therefore follows directly from the number of LLE segments:
\begin{equation}
    N_{\mathrm{TCM}}=1+2N_p+2=2N_{\mathrm{seg}}+1.
    \label{eq:maneuver_count}
\end{equation}
The four-, five-, and eight-segment transfers in Fig~\ref{fig.LLE_spatial} therefore require 9, 11, and 17 planned correction opportunities. This step removes the discrete variable $N_{\mathrm{TCM}}$ from the stochastic search.

\FloatBarrier

\subsection{Reduced Correction-Epoch Parameterization}\label{subsec4.2}
Optimizing all \(2N_{\mathrm{seg}}+1\) correction epochs independently would result in a high-dimensional stochastic optimization problem, with substantial computational cost and increased difficulty in achieving convergence. The segmentation in Section~\ref{subsec3.3} provides the required reduction. Each segment contains only one downstream sensitivity event and can be treated as a local direct-transfer problem. Its error map is still nonlinear, but it is less affected by repeated stretching than a map propagated across the complete multi-loop transfer.

This local structure allows the spacing-ratio idea of Breakwell \cite{breakwell1962spacing,ni2016breakwell} to be used as an  parameterization within each segment. Breakwell generates successive epochs from a geometric ratio of the remaining time to one final target. Here, the same idea is applied separately within every LLE-guided segment. Only the first epoch is optimized, while the later epoch is generated from its remaining time to the next segment boundary. One ratio $c_i$ is therefore assigned to $\mathcal S_i=[t_{\mathrm s}^{(i)},t_{\mathrm e}^{(i)}]$:
\begin{equation}
\begin{aligned}
t_1^{(i)} &= t_{\mathrm{s}}^{(i)}
 +c_i\left(t_{\mathrm{e}}^{(i)}-t_{\mathrm{s}}^{(i)}\right),\\
t_2^{(i)} &= t_{\mathrm{e}}^{(i)}
 -\frac{t_{\mathrm{e}}^{(i)}-t_1^{(i)}}{\rho},
 \qquad i=1,\ldots,N_{\mathrm{seg}}.
\end{aligned}
\label{eq:regular_schedule}
\end{equation}
Both corrections target $t_{\mathrm e}^{(i)}$. The ratio $\rho=3$ gives
\begin{equation}
    t_{\mathrm e}^{(i)}-t_q^{(i)}
    =\frac{t_{\mathrm e}^{(i)}-t_1^{(i)}}
    {\rho^{q-1}},
    \qquad q=1,2,
    \label{eq:remaining_time_rule}
\end{equation}
so the second correction retains one third of the first correction's remaining propagation time. The same recurrence generates the third correction in the terminal segment:
\begin{equation}
t_3^{(N_{\mathrm{seg}})}=t_{\mathrm{e}}^{(N_{\mathrm{seg}})}
-\frac{t_{\mathrm{e}}^{(N_{\mathrm{seg}})}-t_2^{(N_{\mathrm{seg}})}}{\rho}.
\label{eq:final_schedule}
\end{equation}

With this parameterization, only one timing parameter \(c_i\) is optimized for each segment, while the remaining correction epochs are generated automatically from the prescribed spacing ratio. The complete correction schedule is therefore determined by the vector \(\mathbf{c}=[c_1,\ldots,c_{N_{\mathrm{seg}}}]^{\mathrm T}\). This reduces the number of optimization variables from \(2N_{\mathrm{seg}}+1\) correction epochs to \(N_{\mathrm{seg}}\) segment-wise timing parameters, which substantially reduces the dimension of the stochastic optimization problem. For T1--T6, the resulting search dimensions are 4, 5, and 8, instead of 9, 11, and 17, respectively.

\subsection{Second-Order Correction at a Prescribed Epoch}\label{sec:nonlinear_targeting}
Once $\mathbf c$ is specified, the correction law determines the impulse at every generated epoch. Let $\delta\mathbf X_k$ be the navigation-state deviation immediately before correction $k$. An impulsive command changes only the velocity components:
\begin{equation}
\delta\mathbf X_k^{+}=\delta\mathbf X_k+\mathbf B_v\Delta\mathbf v_k,
\qquad
\mathbf B_v=\begin{bmatrix}\mathbf 0_{2\times2}\\\mathbf I_{2\times2}\end{bmatrix}.
\label{eq:post_maneuver_state}
\end{equation}

The state deviation at the target is expanded about the nominal trajectory. Retaining the first- and second-order state-transition terms gives \cite{boone2021orbital,yin2023midcourse}
\begin{equation}
\delta X_{\mathrm{tar},a}
\simeq
\sum_{b=1}^{4}\Phi_{ab}\delta X_{k,b}^{+}
+\frac{1}{2}\sum_{b=1}^{4}\sum_{c=1}^{4}
\Psi_{abc}\delta X_{k,b}^{+}\delta X_{k,c}^{+},
\qquad a=1,\ldots,4.
\label{eq:second_order_state_map}
\end{equation}
Here, $\mathbf\Phi$ and $\mathbf\Psi$ are integrated from $t_k$ to the downstream target using Eq.~\eqref{eq:stm_stt_dynamics}. Each correction has two in-plane control components and imposes two terminal-position conditions. With $\mathbf P_r=[\mathbf I_{2\times2},\mathbf 0_{2\times2}]$, the residual becomes
\begin{equation}
\mathbf e_k(\Delta\mathbf v_k)
=\mathbf P_r\left(
\mathbf\Phi\delta\mathbf X_k^{+}
+\frac{1}{2}\mathbf\Psi
[\delta\mathbf X_k^{+},\delta\mathbf X_k^{+}]
\right),
\label{eq:target_residual}
\end{equation}
where $\mathbf\Psi[\mathbf u,\mathbf v]$ denotes the vector obtained by contracting the last two tensor indices with $\mathbf u$ and $\mathbf v$. The first-order impulse provides the initial iterate,
\begin{equation}
\Delta\mathbf v_k^{(0)}
=-\left(\mathbf P_r\mathbf\Phi\mathbf B_v\right)^{\dagger}
\mathbf P_r\mathbf\Phi\delta\mathbf X_k,
\label{eq:first_order_initial_impulse}
\end{equation}
where $(\cdot)^{\dagger}$ denotes a regularized inverse when the sensitivity matrix is poorly conditioned.

The second-order residual is then minimized by a damped Levenberg--Marquardt (LM) iteration. Its Jacobian is
\begin{equation}
\mathbf J_k
=\mathbf P_r\left(
\mathbf\Phi\mathbf B_v
+\mathbf\Psi[\delta\mathbf X_k^{+},\mathbf B_v]
\right),
\label{eq:target_jacobian}
\end{equation}
and the command update is
\begin{equation}
\delta(\Delta\mathbf v_k)
=-(\mathbf J_k^{\mathrm T}\mathbf J_k+\lambda\mathbf I)^{-1}
\mathbf J_k^{\mathrm T}\mathbf e_k.
\label{eq:lm_update}
\end{equation}
The damping decreases when the trial step reduces the target residual and increases otherwise. The iteration stops when the position residual or the command increment reaches its prescribed tolerance. In each stochastic realization, navigation error enters $\delta\mathbf X_k$, while magnitude and pointing errors are applied to the converged command before the next propagation interval. The command is limited to \SI{30}{\metre\per\second}. The resulting STT--LM targeting procedure for a prescribed correction epoch and downstream target is summarized in Algorithm~\ref{alg:stt_guidance}.

\begin{algorithm}[htbp]
\LinesNumbered
\caption{Second-order STT--LM targeting}
\label{alg:stt_guidance}
\KwIn{Correction epoch $t_k$, target epoch $t_{\mathrm{tar}}$, navigation state $\mathbf{S}_{\mathrm{nav}}$, nominal state $\mathbf{S}_{\mathrm{nom}}(t_k)$, tolerance $\epsilon_{\mathrm{LM}}$, damping $\lambda$, and scale $\nu$}
Set $\delta\mathbf{x}_0=\mathbf{S}_{\mathrm{nav}}-\mathbf{S}_{\mathrm{nom}}(t_k)$\;
Integrate $\mathbf{\Phi}$ and $\mathbf{\Psi}$ from $t_k$ to $t_{\mathrm{tar}}$ using Eq.~\eqref{eq:stm_stt_dynamics}\;
Extract the position-level first- and second-order mappings\;
Initialize $\Delta\mathbf{v}$ using Eq.~\eqref{eq:first_order_initial_impulse}\;
Evaluate the second-order terminal position residual $\mathbf{e}_f$\;
\While{$\|\mathbf{e}_f\|>\epsilon_{\mathrm{LM}}$}{
    Form the residual Jacobian from Eq.~\eqref{eq:target_jacobian}\;
    Compute the LM update from Eq.~\eqref{eq:lm_update}\;
    Evaluate the trial residual after $\Delta\mathbf{v}_{\mathrm{try}}=\Delta\mathbf{v}+\delta(\Delta\mathbf{v})$\;
    \eIf{the trial residual norm decreases}{
        Accept the trial impulse and set $\lambda\leftarrow\lambda/\nu$\;
    }{
        Reject the trial impulse and set $\lambda\leftarrow\lambda\nu^2$\;
    }
}
\Return{$\Delta\mathbf{v}$}\;
\end{algorithm}

\subsection{Chance-Constrained Correction Schedule Optimization}\label{subsec4.4}
After the maneuver architecture and epoch parameterization are fixed, the correction strategy is completely determined by the reduced timing vector \(\mathbf{c}\). Its performance under uncertainty must be evaluated through nonlinear trajectory propagation, so neither the expected correction cost nor the mission-success probability is available in closed form. A sample-average approximation \cite{calafiore2006scenario,campi2008exact} is therefore adopted to formulate the chance-constrained schedule optimization. 

A fixed set of standardized disturbance samples, \(\Xi_M=\{\boldsymbol{\xi}_m\}_{m=1}^{M}\), is generated before optimization. Each realization contains four injection-error components and, for every correction, four navigation-error and two maneuver-execution-error components, giving a total dimension of \(4+6N_{\mathrm{TCM}}\). The standardized samples are transformed according to the prescribed uncertainty model during trajectory propagation. The same disturbance set \(\Xi_M\) is used for all candidate schedules. This common-random-number strategy reduces sampling noise in candidate comparisons, so differences in the estimated objective and reliability primarily reflect changes in the correction schedule rather than changes in the sampled disturbances. For each realization \(m\), let \(I_m(\mathbf c)\) denote the joint-success indicator, with \(I_m=1\) if the joint-success event defined in Section~\ref{subsec2.3} is satisfied and \(I_m=0\) otherwise. The specific success criteria and their numerical thresholds are given in Section~\ref{subsec5.1}. The sample-based estimate of the mission-success probability is
\begin{equation}
\widehat p_M(\mathbf{c})=\frac{1}{M}\sum_{m=1}^{M}I_m(\mathbf{c}).
\label{eq:sample_reliability}
\end{equation}

Let $\mathcal S(\mathbf c)=\{m\mid I_m(\mathbf c)=1\}$ denote the set of successful realizations. The correction-cost objective is evaluated over this set as:
\begin{equation}
\widehat J_M(\mathbf{c})=
\overline D_{\mathcal{S}}(\mathbf{c})
+3s_{D,\mathcal{S}}(\mathbf{c}),
\label{eq:sample_objective}
\end{equation}
where $\overline D_{\mathcal{S}}$ and $s_{D,\mathcal{S}}$ are the sample mean and standard deviation over successful realizations. Equation~\eqref{eq:general_stochastic_problem} is therefore specialized as
\begin{equation}
\begin{aligned}
\min_{\mathbf{c}}\quad &\widehat J_M(\mathbf{c}),\\
\mathrm{s.t.}\quad
&\widehat p_M(\mathbf{c})\geq p_{\mathrm{req}},
\qquad p_{\mathrm{req}}=0.95,\\
&0.01\leq c_i\leq0.50,
\qquad i=1,\ldots,N_{\mathrm{seg}},\\
&t_{k+1}(\mathbf{c})>t_k(\mathbf{c}),\qquad
t_{\mathrm{tar},k}(\mathbf{c})>t_k(\mathbf{c}).
\end{aligned}
\label{eq:sample_so}
\end{equation}
The chance constraint directly enforces the joint mission-success requirement rather than embedding individual failure criteria in a weighted penalty function. The correction-cost statistics are computed only from successful realizations, while all \(M\) samples, including failed realizations, contribute to the reliability estimate. In this way, failed trajectories do not introduce artificial correction-cost values into the objective, but their effect on mission reliability is fully retained through the chance constraint.

The resulting stochastic optimization problem is solved using a two-stage surrogate-based search. A scrambled Sobol sequence is transformed to standard normal variables and used to form $M_{\mathrm g}=256$ global-search scenarios \cite{joe2008sobol}. These samples form a subset of a larger \(M_r=512\) set used for local refinement, so the same disturbances are retained as the sample size increases. The initial design contains the current strategy, positive and negative perturbations of each timing coordinate, and 64 Sobol points in the bounded design space. A constrained radial-basis-function (RBF) surrogate then explores the space for at most 600 candidate evaluations \cite{regis2007rbf}. At each evaluation, the correction epochs are generated from Eqs.~\eqref{eq:regular_schedule} and \eqref{eq:final_schedule}, and the associated stochastic performance is evaluated using the STT--LM targeting procedure.

The best global candidate is re-evaluated on the 512-scenario set. Pattern search then performs a local constrained refinement with a budget of 300 evaluations \cite{audet2006mads}. The initial, global, and locally refined candidates are compared on this common set. Feasibility is checked first. Among candidates satisfying the reliability and timing constraints, the one with the smallest $\widehat J_M$ is retained. If no candidate is feasible, the most reliable candidate is archived for diagnosis and the existing strategy is preserved.

To avoid accepting a strategy that performs well only on the optimization samples, the retained candidate is further evaluated using an independent validation set of \(M_v=500\) realizations. For \(n_s\) successful realizations, the one-sided 95\% Clopper--Pearson lower confidence bound is \cite{clopper1934interval}
\begin{equation}
p_{\mathrm L}=F_{\mathrm{Beta}}^{-1}
\left(0.05;\,n_{\mathrm s},M_{\mathrm v}-n_{\mathrm s}+1\right).
\label{eq:binomial_lower}
\end{equation}
The candidate and the current strategy are evaluated with the same 500 validation disturbances. The candidate is accepted only when $p_{\mathrm L}\geq0.95$, all epochs are operationally constraints, and its validation success rate is not lower than that of the current strategy. Otherwise, the candidate is rejected and the existing timing schedule is retained. This independent validation step separates strategy optimization from the statistical evidence used for reliability acceptance. The complete optimization and validation procedure is summarized in Algorithm~\ref{alg:robust_optimization}.

\begin{algorithm}[htbp]
\LinesNumbered
\caption{Chance-constrained optimization and independent acceptance}
\label{alg:robust_optimization}
\KwIn{Nominal trajectory, LLE boundaries, uncertainty model, success criteria, timing bounds, and sample sizes $M_{\mathrm g}$, $M_{\mathrm r}$, and $M_{\mathrm v}$}
\KwOut{Accepted timing vector $\mathbf{c}^{*}$, correction sequence $\mathcal{T}^{*}$, and validation diagnostics}
Construct $2N_{\mathrm{seg}}+1$ epochs from Eqs.~\eqref{eq:regular_schedule} and \eqref{eq:final_schedule}\;
Generate nested Sobol sets $\Xi_{\mathrm g}\subset\Xi_{\mathrm r}$\;
Build the initial candidates from the existing strategy, coordinate perturbations, and scrambled Sobol points\;
\ForEach{candidate $\mathbf{c}$ proposed by the constrained RBF search}{
  Generate $\mathcal{T}(\mathbf{c})$ and the target boundary of each correction\;
  \ForEach{$\boldsymbol{\xi}_m\in\Xi_{\mathrm g}$ in parallel}{
    Inject the prescribed errors and apply STT--LM targeting at every epoch\;
    Compute correction cost, terminal errors, path metrics, and maximum pulse\;
    Evaluate the joint-success indicator $I_m(\mathbf{c})$\;
  }
  Update $\widehat p_M$, $\widehat J_M$, and the RBF surrogate\;
}
Re-evaluate the global candidate on $\Xi_{\mathrm r}$ and refine it by pattern search\;
Select $\mathbf{c}^{*}$ using the feasibility-first rule\;
Evaluate $\mathbf{c}^{*}$ and the existing strategy on independent $\Xi_{\mathrm v}$\;
Compute $p_{\mathrm L}$ from Eq.~\eqref{eq:binomial_lower}\;
\eIf{$p_{\mathrm L}\geq0.95$, timing is feasible, and validation success is not worse}{
  Accept $\mathbf{c}^{*}$ and adopt $\mathcal{T}^{*}$\;
}{
  Archive the candidate and retain the existing timing sequence\;
}
\end{algorithm}

\begin{figure}[htbp]
    \centering
    \includegraphics[width=0.86\textwidth]{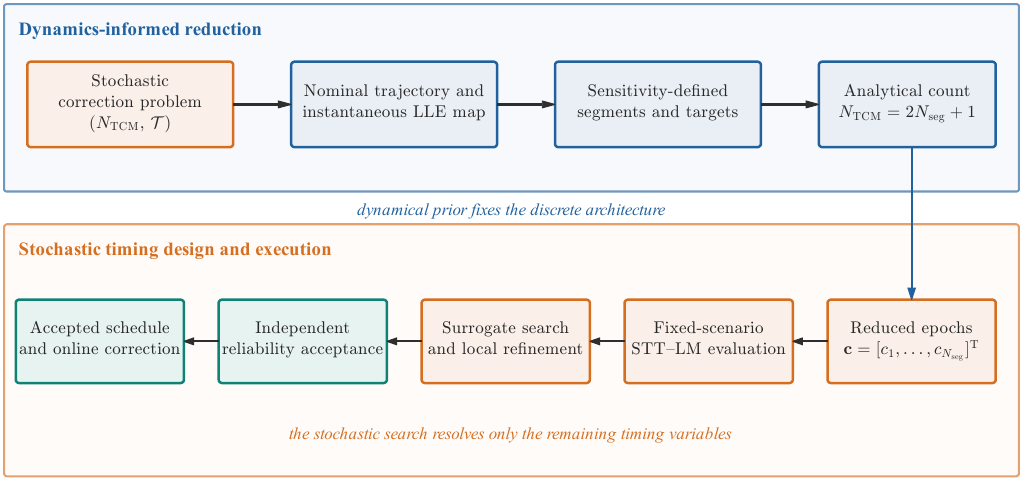}
    \caption{Overall workflow of the proposed LLE-guided segmented mid-course correction strategy.}
    \label{fig:lle_segmented_framework}
\end{figure}

The complete design chain of mid-course correction for low-energy transfers is summarized in Fig.~\ref{fig:lle_segmented_framework}. LLE segmentation resolves the discrete architecture, the stochastic search handles only $N_{\mathrm{seg}}$ timing variables, and independent samples determine whether the resulting schedule is accepted.

\FloatBarrier

\section{Numerical Results and Discussion}\label{sec5}
The proposed LLE-guided segmented correction strategy is evaluated on the six representative low-energy transfers introduced in Section \ref{subsec2.2} and compared with the classical Breakwell spacing-ratio method. The comparison focuses on two questions: whether the proposed strategy improves the probability of recovering a disturbed low-energy transfer, and how the two methods differ in correction cost, terminal accuracy, and preservation of the nominal transfer path. T3 is first used to illustrate the correction performance in detail, followed by paired Monte Carlo (MC) comparisons across all six trajectories.

\subsection{Simulation Setup and Performance Evaluation Framework}\label{subsec5.1}
All trajectories are propagated in the planar PCRTBP with relative and absolute integration tolerances of \(10^{-12}\), and the nonlinear targeting tolerance $\epsilon_{\mathrm{LM}}$ is \(10^{-9}\) in normalized position units. The common uncertainty model is summarized in Table~\ref{tab:error_model}. Injection, navigation, and maneuver-execution errors are modeled as zero-mean Gaussian disturbances. Navigation accuracy switches between the near-Earth and far-field models at a geocentric distance of \(12\,000\) km~\cite{he2007lowenergy}, while independent magnitude and in-plane direction errors are applied to each executed maneuver.

\begin{table}[htbp]
	\caption{Uncertainty model for mid-course correction design and paired evaluation.}
	\label{tab:error_model}
	\centering
	\small
	\begin{tabular}{lccc}
		\toprule
		Error source & Position (3$\sigma$) (km) & Velocity (3$\sigma$) (m/s) & Execution dispersion (3$\sigma$)\\
		\midrule
		Injection & 3.000 & 10.0 & --\\
		Near-Earth navigation & 0.060 & 0.3 & --\\
		Far-field navigation & 0.300 & 0.3 & --\\
		Maneuver execution & -- & -- & 1\% $\|\Delta \mathbf v\|$; $0.2^{\circ}$ direction\\
		\bottomrule
	\end{tabular}
\end{table}

The reference strategy follows Breakwell's remaining-time spacing rule. The maneuver count is estimated from the initially projected terminal position error \(\delta\mathbf r_{f,0}\) and the required terminal accuracy \(\delta\mathbf r_f^\ast\) as

\begin{equation}
	N_{\mathrm B}\approx 1+\log_3\left(
	\frac{\|\delta\mathbf{r}_{f,0}\|}{\|\delta\mathbf{r}_f^{*}\|}
	\right).
	\label{eq:breakwell_count}
\end{equation}
the correction epochs satisfy
\begin{equation}
	\frac{\tau_1}{\tau_2}=\frac{\tau_2}{\tau_3}=\cdots=
	\frac{\tau_{N_{\mathrm B}-1}}{\tau_{N_{\mathrm B}}}=3.
	\label{eq:breakwell_spacing}
\end{equation}
Here, $\tau_k=t_f-t_{k,\mathrm B}$ denotes the remaining flight time at the $k$th Breakwell correction. The first correction epoch of the Breakwell method is set to be the same as that of the accepted Segmented schedule. Subsequent Breakwell epochs satisfy \(t_{k+1,\mathrm B}=t_f-(t_f-t_{k,\mathrm B})/3\). Both methods use the same uncertainty model, maneuver limits, numerical tolerances, and STT--LM targeting law.

For each trajectory, the two correction strategies are evaluated using 1000 paired Monte Carlo realizations. Within each pair, both methods are subjected to the same standardized injection, navigation, and indexed maneuver-execution errors, so their performance is compared under identical disturbances. The samples used for schedule optimization, independent reliability validation, and final performance evaluation are mutually disjoint. For a direct transfer, correction performance can generally be characterized by the terminal position and velocity errors. For a low-energy transfer, however, the nominal path is itself an important part of the mission design because it exploits multi-body dynamics to achieve the intended transfer characteristics and satisfy trajectory constraints. Consequently, a corrected trajectory may meet the terminal arrival requirements while deviating substantially from the designed low-energy path. Terminal errors alone are therefore insufficient to evaluate whether the original transfer characteristics are preserved. The time-domain root-mean-square position deviation is then introduced to quantify the overall deviation from the nominal trajectory as

\begin{equation}
	G_{\text{path}}
	=\sqrt{\frac{1}{t_f-t_0}\int_{t_0}^{t_f}
		\left\|\mathbf{r}(t)-\mathbf{r}^*(t)\right\|^2\,dt},
	\label{eq:path_rmse}
\end{equation}
and the corresponding path-recovery rate is defined as
\begin{equation}
	\eta_{\text{path}}
	=\frac{G_{\text{path}}^{\text{uncorr}}-G_{\text{path}}^{\text{corr}}}
	{G_{\text{path}}^{\text{uncorr}}}.
	\label{eq:path_improvement}
\end{equation}
The joint-success event $\mathcal{E}$ requires all terminal, correction-cost, pulse, full-path RMSE, and path-recovery limits to be satisfied, together with successful propagation and targeting. It thus tests preservation of the intermediate transfer path rather than considering arrival accuracy alone. Detailed joint mission success is defined by the criteria in Table~\ref{tab:success_protocol}. Reliability is evaluated over all 1000 paired realizations, while correction cost, terminal errors, and path metrics are compared only for common-success pairs. This distinction separates the ability to recover a disturbed trajectory from the conditional performance once both methods have succeeded.

\begin{table}[htbp]
	\caption{Joint success criteria used in the paired Monte Carlo evaluation.}
	\label{tab:success_protocol}
	\centering
	\small
	\begin{tabularx}{\textwidth}{>{\raggedright\arraybackslash}p{0.31\textwidth}>{\raggedright\arraybackslash}p{0.25\textwidth}X}
		\toprule
		Criterion & Required value & Physical role\\
		\midrule
		Total correction & $\Delta \mathbf v_{\mathrm{tot}}\leq\SI{300}{\metre\per\second}$ & Mission-level correction budget\\
		Terminal position error & $e_r\leq\SI{50}{\kilo\metre}$ & Arrival-position accuracy\\
		Terminal velocity error & $e_v\leq\SI{50}{\metre\per\second}$ & Arrival-velocity accuracy\\
		Corrected path RMSE & $G_{\mathrm{path}}^{\mathrm{corr}}\leq\SI{100000}{\kilo\metre}$ & Exclusion of gross path divergence\\
		Path improvement & $\eta_{\mathrm{path}}\geq0.95$ & Recovery of at least 95\% of the uncorrected deviation\\
		Individual maneuver & $\max_k\|\Delta\mathbf v_k\|\leq\SI{30}{\metre\per\second}$ & Single-command capability\\
		\bottomrule
	\end{tabularx}
\end{table}

The corresponding accepted relative timing vectors obtained using the simulation settings described above are presented in Table~\ref{tab:accepted_timing_vectors}. Each component specifies the first correction in one LLE segment. The remaining epochs follow Eqs.~\eqref{eq:regular_schedule} and \eqref{eq:final_schedule}.

\begin{table}[htbp]
	\caption{Accepted relative correction-epoch vectors $\mathbf c^{*}$.}
	\label{tab:accepted_timing_vectors}
	\centering
	\small
	\begin{tabular}{cl}
		\toprule
		Case & $\mathbf c^{*}$\\
		\midrule
		T1 & $[0.0102,\,0.0137,\,0.1957,\,0.1471]^{\mathrm T}$\\
		T2 & $[0.0103,\,0.0142,\,0.1932,\,0.1091]^{\mathrm T}$\\
		T3 & $[0.0100,\,0.0838,\,0.0100,\,0.1854,\,0.1206]^{\mathrm T}$\\
		T4 & $[0.0100,\,0.0761,\,0.0100,\,0.1830,\,0.1572]^{\mathrm T}$\\
		T5 & $[0.0149,\,0.0163,\,0.2018,\,0.1532,\,0.1331,\,0.0292,\,0.1662,\,0.1233]^{\mathrm T}$\\
		T6 & $[0.0138,\,0.0473,\,0.2014,\,0.1499,\,0.1261,\,0.0400,\,0.1667,\,0.1595]^{\mathrm T}$\\
		\bottomrule
	\end{tabular}
\end{table}

\subsection{Correction Performance on a Representative Low-Energy Transfer}\label{subsec5.2}

T3 is selected as a representative case because it contains both Earth-centered perigee-to-perigee return arcs and a final Earth-Moon transfer arc, covering the two recurring segment geometries identified in Section~\ref{subsec3.3}. The transfer lasts \SI{40.77}{day} and contains five LLE-guided segments and 11 planned correction opportunities, with segment boundaries at 0, 8.846, 12.734, 21.101, 33.475, and 40.766 days.

\begin{figure}[htbp]
    \centering
    \includegraphics[width=1\textwidth]{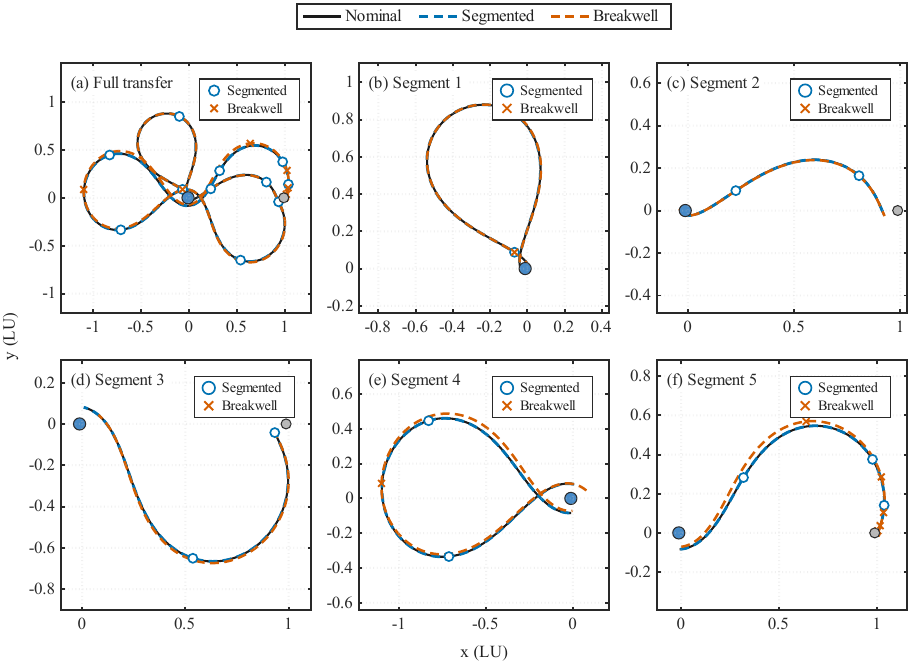}
    \caption{Representative segment-wise comparison of trajectory correction between the Segmented and Breakwell strategies. (a) Complete transfer; (b)--(f) five LLE-guided segments.}
    \label{fig:t3_geometry}
\end{figure}

A representative Monte Carlo realization first illustrates how the two strategies differ after the initial injection error has been removed. During the first two segments, both corrected trajectories remain close to the nominal path because the initial deviation is still effectively controlled. The difference becomes pronounced in the subsequent Earth-centered loops. The Breakwell trajectory gradually departs from the nominal path before the long return arc, and this deviation persists into the final lunar-approach segment. The Segmented strategy, in contrast, repeatedly reduces the accumulated deviation at the intervening sensitivity boundaries and remains close to the nominal geometry throughout the transfer. By the final segment, both strategies satisfy the terminal requirements, but only the Segmented trajectory has recovered the nominal transfer path over the complete  segment (Fig.~\ref{fig:t3_geometry}).

\begin{figure}[htbp]
    \centering
    \begin{subfigure}[t]{0.49\textwidth}
        \centering
        \includegraphics[height=5.6cm]{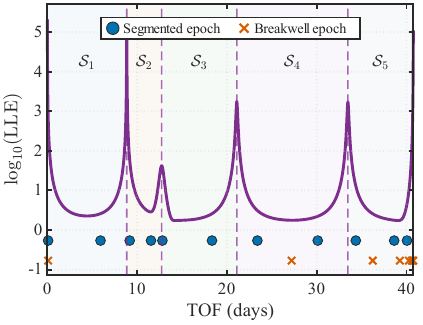}
        \caption{Correction epochs.}
    \end{subfigure}
    \hfill
    \begin{subfigure}[t]{0.49\textwidth}
        \centering
        \includegraphics[height=5.6cm]{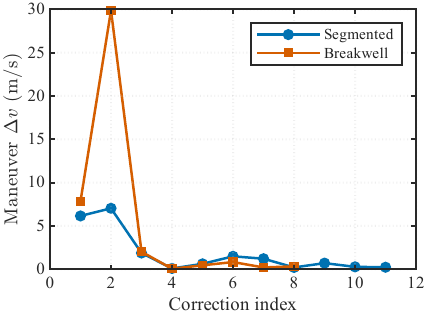}
        \caption{Maneuver magnitudes.}
    \end{subfigure}
    \caption{Representative correction timing and maneuver allocation. (a) Shaded bands denote the five LLE-guided segments, with the correction epochs shown below. (b) Executed maneuver magnitudes in MC sample 247.}
    \label{fig:t3_mechanism}
\end{figure}

Although both strategies apply their first correction at \SI{2.1231}{h} after injection, their subsequent schedules differ markedly (Fig.~\ref{fig:t3_mechanism}). The Segmented strategy distributes its 11 corrections across the five sensitivity-defined intervals, with maneuvers placed before and after the internal segment boundaries. The upstream corrections reduce the errors entering the sensitive regions, while the downstream corrections reduce the residual deviations carried into the following segments. In contrast, the eight Breakwell corrections follow a global remaining-time spacing rule and become increasingly concentrated near arrival, without explicitly accounting for the intermediate sensitivity peaks. Errors that are not sufficiently reduced before these peaks can undergo further amplification, increasing the demand for subsequent trajectory recovery. The segmented architecture instead limits this accumulation through successive local corrections.

For this realization, the improved path recovery is achieved with a lower total correction cost despite the larger number of maneuvers. The Segmented strategy requires \SI{19.85}{\metre\per\second}, compared with \SI{41.59}{\metre\per\second} for Breakwell. The corresponding full-path position RMSE values are \SI{246.56}{\kilo\metre} and \SI{7038.44}{\kilo\metre}, respectively, while the uncorrected full-path RMSE is \SI{171642.20}{\kilo\metre}. Both strategies satisfy the terminal constraints: the terminal position errors are \SI{14.22}{\kilo\metre} for Segmented and \SI{0.21}{\kilo\metre} for Breakwell, while the terminal velocity errors are \SI{9.09}{\metre\per\second} and \SI{20.94}{\metre\per\second}, respectively. Thus, Breakwell achieves a smaller terminal position error but incurs a substantially larger intermediate deviation and a higher correction cost. The Segmented strategy, by contrast, maintains the required arrival accuracy while recovering the nominal path more effectively. These observations illustrate the correction performance for a single realization. The statistical comparison is presented in Section~\ref{subsec5.3}.

\begin{figure}[!ht]
	\centering
	\includegraphics[width=0.8\textwidth]{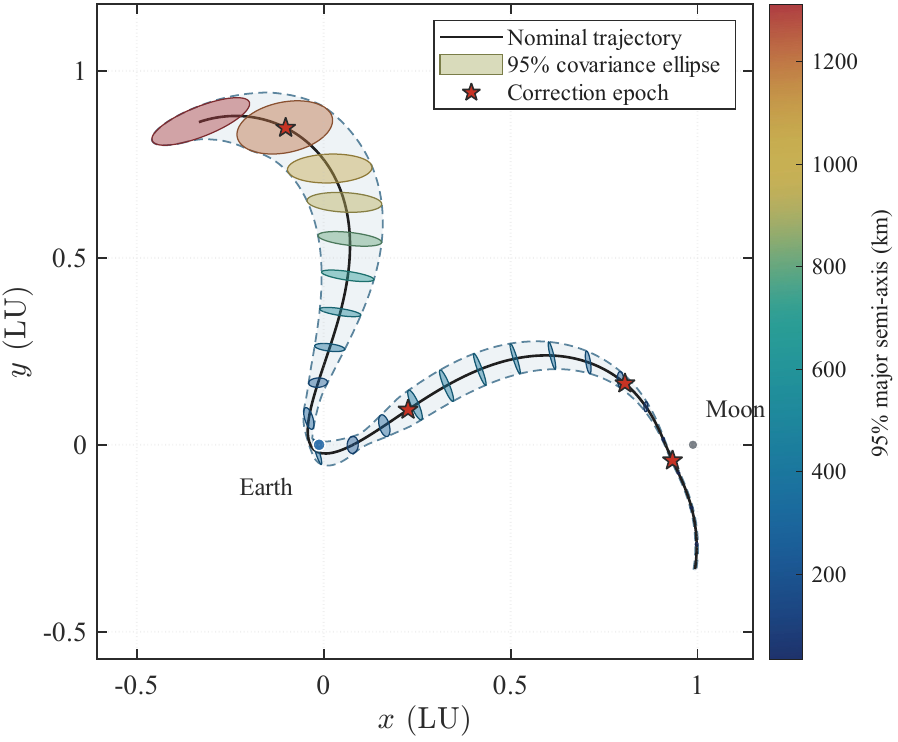}
	\caption{Position-error evolution under the Segmented correction strategy.}
	\label{fig:t3_error_convergence}
\end{figure}

The covariance evolution provides an ensemble-level view of the same mechanism. As the trajectory approaches a high-sensitivity region, the position-error ellipses expand as the accumulated uncertainty is amplified. The subsequent correction changes the downstream state evolution, after which the dispersion contracts and the trajectory is brought back toward the nominal path. This repeated expansion--correction--contraction behavior is visible across successive segments in Fig.~\ref{fig:t3_error_convergence}. The segmented architecture therefore acts by repeatedly limiting the error before it can be carried through multiple sensitive passages, rather than attempting to recover the accumulated deviation near the end of the transfer.

\FloatBarrier
\subsection{Paired Cross-Trajectory Performance}\label{subsec5.3}
The T3 example explains how sensitivity-guided correction can suppress the accumulation of trajectory error. The remaining question is whether this mechanism produces a consistent advantage across different low-energy transfer geometries. The paired Monte Carlo results for T1--T6 are examined below from three complementary perspectives: recovery reliability, the trade-off between terminal accuracy and full-path preservation, and correction demand.

\subsubsection{Reliability and Paired Success Outcomes}\label{subsubsec5.3.1}
Before the final Monte Carlo comparison, all optimized schedules are required to pass the independent reliability-acceptance test introduced in Section~4.4. The four-, five-, and eight-segment transfers correspond to 9, 11, and 17 planned corrections, respectively. Across all six cases, 491--498 of the 500 independent validation samples are successful, and the smallest one-sided 95\% Clopper--Pearson lower confidence bound is 96.88\%. The accepted architectures and validation results are summarized in Table~\ref{tab:optimized_schedules}, confirming that each schedule satisfies the prescribed reliability requirement before entering the final test.

\begin{table}[!ht]
  \caption{Correction architecture and independent schedule acceptance.}
  \label{tab:optimized_schedules}
  \centering
  \small
  \begin{tabular}{cccccc}
    \toprule
    & & \multicolumn{2}{c}{Correction architecture} & \multicolumn{2}{c}{Independent acceptance}\\
    \cmidrule(lr){3-4}\cmidrule(lr){5-6}
    Case & TOF (day) & $N_{\mathrm{seg}}$ & $N_{\mathrm{TCM}}$ & Successes & $p_{\mathrm L}$\\
    \midrule
    T1 & 31.92 & 4 & 9  & 498/500 & 98.75\%\\
    T2 & 31.95 & 4 & 9  & 494/500 & 97.65\%\\
    T3 & 40.77 & 5 & 11 & 491/500 & 96.88\%\\
    T4 & 40.79 & 5 & 11 & 498/500 & 98.75\%\\
    T5 & 58.45 & 8 & 17 & 491/500 & 96.88\%\\
    T6 & 58.53 & 8 & 17 & 491/500 & 96.88\%\\
    \bottomrule
  \end{tabular}
\end{table}

\FloatBarrier

\begin{figure}[htbp]
    \centering
    \begin{subfigure}[t]{0.6\textwidth}
        \includegraphics[width=\textwidth]{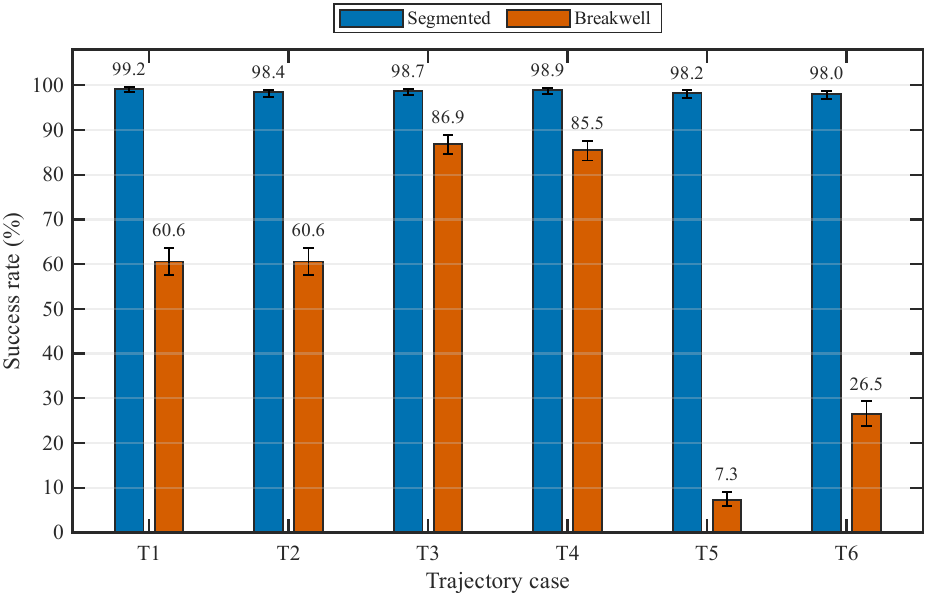}
        \caption{Success probability.}
    \end{subfigure}
    \hfill
    \begin{subfigure}[t]{0.6\textwidth}
        \includegraphics[width=\textwidth]{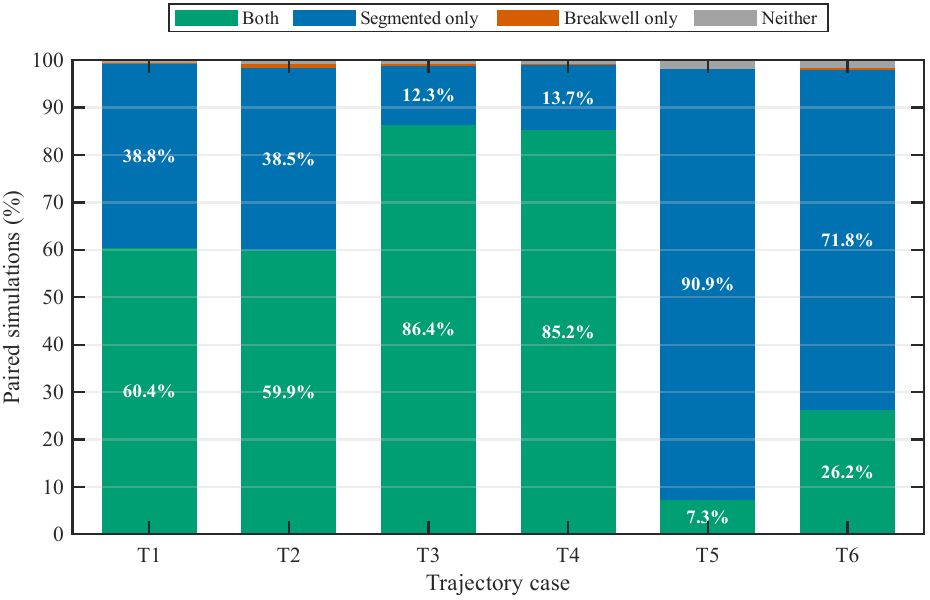}
        \caption{Paired outcomes.}
    \end{subfigure}
    \caption{Reliability over 1000 paired realizations per trajectory. (a) Success rates with Wilson 95\% intervals. (b) Joint outcome counts for the same standardized disturbances.}
    \label{fig:reliability_comparison}
\end{figure}

The principal statistical advantage of the proposed strategy is its ability to recover a much larger fraction of the disturbed transfers. Across T1--T6, the Segmented success rate remains within 98.0--99.2\%, while Breakwell reaches only 60.6--86.9\% on T1--T4, decreases to 7.3\% on T5, and reaches 26.5\% on T6 (Fig.~\ref{fig:reliability_comparison}(a)). The contrast between T5 and T6 is particularly informative. Their flight times and numbers of LLE segments are nearly identical, yet the Breakwell success rates differ by 19.2 percentage points. Correction difficulty therefore cannot be characterized by transfer duration or maneuver count alone; the location and sequence of the sensitive passages are also important.

The paired comparison uses the same standardized disturbance samples for both strategies. Under identical realizations, Segmented-only successes range from 123 for T3 to 909 for T5, while Breakwell-only successes never exceed seven (Fig.~\ref{fig:reliability_comparison}(b)). The strong asymmetry shows that the proposed strategy primarily recovers disturbances that are not successfully corrected by the Breakwell schedule, rather than merely exchanging one set of successful realizations for another.

The failure modes further explain where this difference arises. A path violation is defined by $\eta_{\mathrm{path}}<0.95$, a position violation by $e_r>\SI{50}{\kilo\metre}$, a velocity violation by $e_v>\SI{50}{\metre\per\second}$, and a pulse violation by $\max_k\|\Delta\mathbf v_k\|>\SI{30}{\metre\per\second}$. Because one realization can violate several criteria, these categories are not mutually exclusive. Inadequate path recovery is prevalent among unsuccessful Breakwell realizations, occurring in 93.7--100\% of them and frequently accompanied by terminal-error violations on T5--T6 (Table~\ref{tab:failure_modes}). By comparison, the Segmented strategy produces only 8--20 failures per trajectory, without a common dominant failure pattern. All applied pulses respect the \SI{30}{\metre\per\second} limit.

\begin{table}[htbp]
  \caption{Criterion-violation rates among each method's unsuccessful realizations.}
  \label{tab:failure_modes}
  \centering
  \small
  \setlength{\tabcolsep}{3.5pt}
  \begin{tabular}{clrrrrrrr}
    \toprule
    Case & Method & $n_{\mathrm f}$ & Path (\%) & Position (\%) & Velocity (\%) & Total $\Delta v$ (\%) & Pulse (\%) \\
    \midrule
    T1 & Segmented & 8 & 75.0 & 50.0 & 25.0 & 0.0 & 0.0 \\
     & Breakwell & 394 & 100.0 & 0.5 & 0.5 & 0.0 & 0.0 \\
    \addlinespace[2pt]
    T2 & Segmented & 16 & 37.5 & 75.0 & 12.5 & 0.0 & 0.0 \\
     & Breakwell & 394 & 99.2 & 0.5 & 17.3 & 0.0 & 0.0 \\
    \addlinespace[2pt]
    T3 & Segmented & 13 & 46.2 & 61.5 & 7.7 & 7.7 & 0.0 \\
     & Breakwell & 131 & 100.0 & 5.3 & 5.3 & 0.0 & 0.0 \\
    \addlinespace[2pt]
    T4 & Segmented & 11 & 63.6 & 45.5 & 9.1 & 9.1 & 0.0 \\
     & Breakwell & 145 & 100.0 & 6.2 & 6.2 & 0.0 & 0.0 \\
    \addlinespace[2pt]
    T5 & Segmented & 18 & 66.7 & 33.3 & 0.0 & 0.0 & 0.0 \\
     & Breakwell & 927 & 93.7 & 99.5 & 99.6 & 0.0 & 0.0 \\
    \addlinespace[2pt]
    T6 & Segmented & 20 & 55.0 & 45.0 & 0.0 & 0.0 & 0.0 \\
     & Breakwell & 735 & 95.9 & 76.3 & 79.7 & 0.0 & 0.0 \\
    \bottomrule
  \end{tabular}
\end{table}

\FloatBarrier

\subsubsection{Terminal Accuracy and Path Preservation}\label{subsubsec5.3.3}

\begin{figure}[!htbp]
	\centering
	\includegraphics[width=0.84\textwidth]{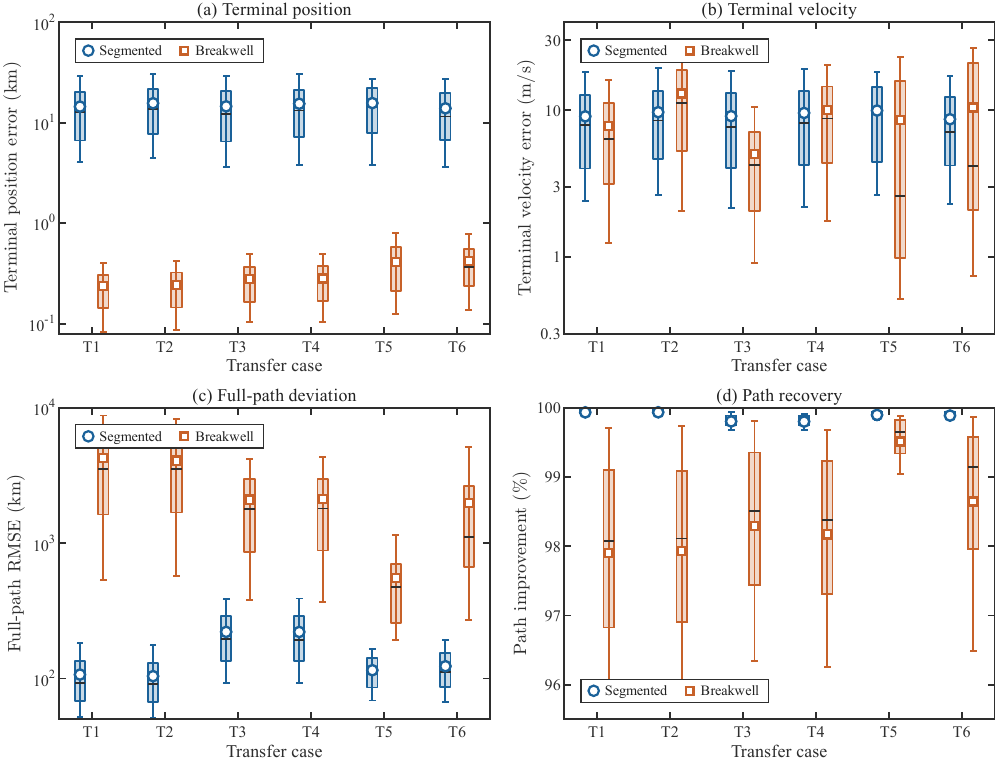}
	\caption{Accuracy and path preservation on common-success pairs. Boxes span the interquartile range, horizontal strokes show medians, open markers show means, and whiskers span the 10th--90th percentiles.}
	\label{fig:accuracy_path}
\end{figure}

The two strategies exhibit a clear trade-off between terminal precision and control of the complete transfer. Among the common-success realizations, Breakwell generally achieves smaller terminal position errors than the Segmented strategy. This behavior follows naturally from its remaining-time spacing rule: successive correction epochs are progressively compressed toward the arrival epoch, so the final state is repeatedly retargeted over increasingly short propagation intervals. The Segmented strategy instead distributes its correction opportunities throughout the transfer according to the sensitivity structure, with the primary objective of preventing large deviations from accumulating through successive high-LLE passages. The different terminal-error distributions in Fig.~\ref{fig:accuracy_path}(a) therefore reflect two different correction philosophies rather than a simple difference in targeting capability.

The terminal accuracy obtained by the Segmented strategy is nevertheless comfortably within the prescribed mission requirement. Its median terminal position error ranges from \SI{11.52}{\kilo\metre} to \SI{14.64}{\kilo\metre}, corresponding to only 23.0--29.3\% of the \SI{50}{\kilo\metre} limit. The terminal velocity errors also remain within the required \SI{50}{\metre\per\second} bound (Fig.~\ref{fig:accuracy_path}(b)). The additional terminal precision obtained by Breakwell therefore lies beyond what is required by the adopted mission-success criterion.

A more consequential difference appears when the complete transfer is considered. The Segmented strategy yields median full-path RMSE values of \SI{90.59}{\kilo\metre}--\SI{196.13}{\kilo\metre}, while Breakwell gives \SI{472.28}{\kilo\metre}--\SI{3552.94}{\kilo\metre} for the same common-success pairs (Fig.~\ref{fig:accuracy_path}(c)). This corresponds to a 76.2--97.4\% reduction. The path-recovery metric leads to the same conclusion: Segmented removes 99.82--99.95\% of the uncorrected path deviation, compared with 98.08--99.65\% for Breakwell (Fig.~\ref{fig:accuracy_path}(d)). Although the normalized percentages appear close, a difference of approximately one percentage point can correspond to hundreds or even thousands of kilometers of full-path deviation on these low-energy transfers. The smaller terminal error of Breakwell is therefore accompanied by substantially poorer control of the intermediate trajectory, while the Segmented strategy preserves the nominal low-energy path while still satisfying the required arrival accuracy.

\subsubsection{Correction Cost and Maneuver Characteristics}\label{subsubsec5.3.2}
\begin{figure}[htbp]
    \centering
    \includegraphics[width=0.92\textwidth]{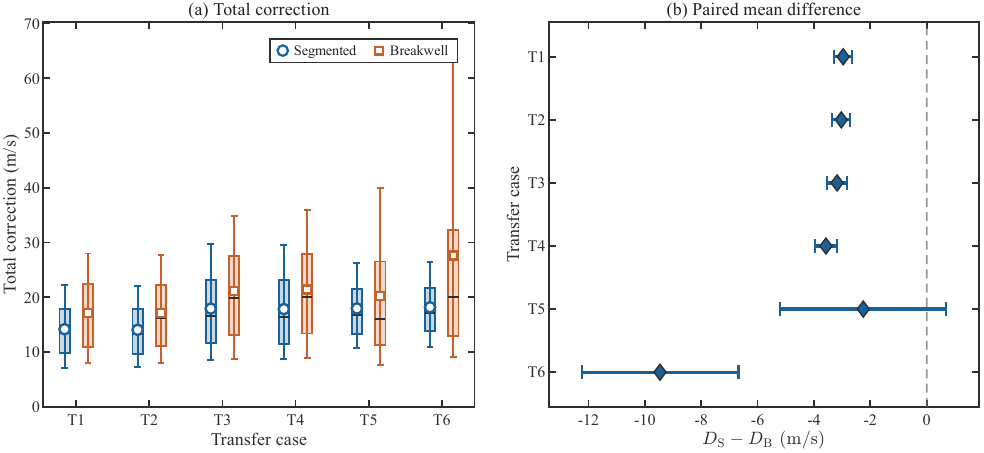}
    \caption{Correction-cost comparison on common-success pairs:
    	(a) total correction distributions; (b) paired differences
    	$D_{\mathrm{S}}-D_{\mathrm{B}}$ with 95\% confidence intervals.}
    \label{fig:correction_cost}
\end{figure}

The higher recovery reliability does not require a larger correction budget. For the common-success pairs, the Segmented strategy has a lower mean correction cost for every transfer, ranging from \SI{14.07}{\metre\per\second} to \SI{18.18}{\metre\per\second}, compared with \SI{17.10}{\metre\per\second}--\SI{27.64}{\metre\per\second} for Breakwell. The resulting mean reduction ranges from 11.13\% for T5 to 34.24\% for T6. The Segmented distributions are also consistently more concentrated, with narrower interquartile and 10th--90th percentile ranges (Fig.~\ref{fig:correction_cost}(a)). The lower dispersion indicates that its correction demand is less sensitive to the particular disturbance realization, giving a more consistent maneuver requirement under uncertainty. The T5 distribution should nevertheless be interpreted together with its reliability result because only 73 common-success pairs remain after the large number of Breakwell failures.

Let $D_{\mathrm S}$ and $D_{\mathrm B}$ denote the total correction costs of the Segmented and Breakwell strategies for the same common-success realization. Their paired differences provide a direct comparison under identical disturbances. The mean value of $D_{\mathrm S}-D_{\mathrm B}$ is negative for all six transfers (Fig.~\ref{fig:correction_cost}(b)), and 70.01\% of the 3254 common-success pairs require less correction with the Segmented strategy. The pooled mean and median differences are \SI{-3.70}{\metre\per\second} and \SI{-2.38}{\metre\per\second}, respectively. Thus, the reduction observed in the marginal distributions also persists on a realization-by-realization basis.

\begin{figure}[!htbp]
	\centering
	\includegraphics[width=0.92\textwidth]{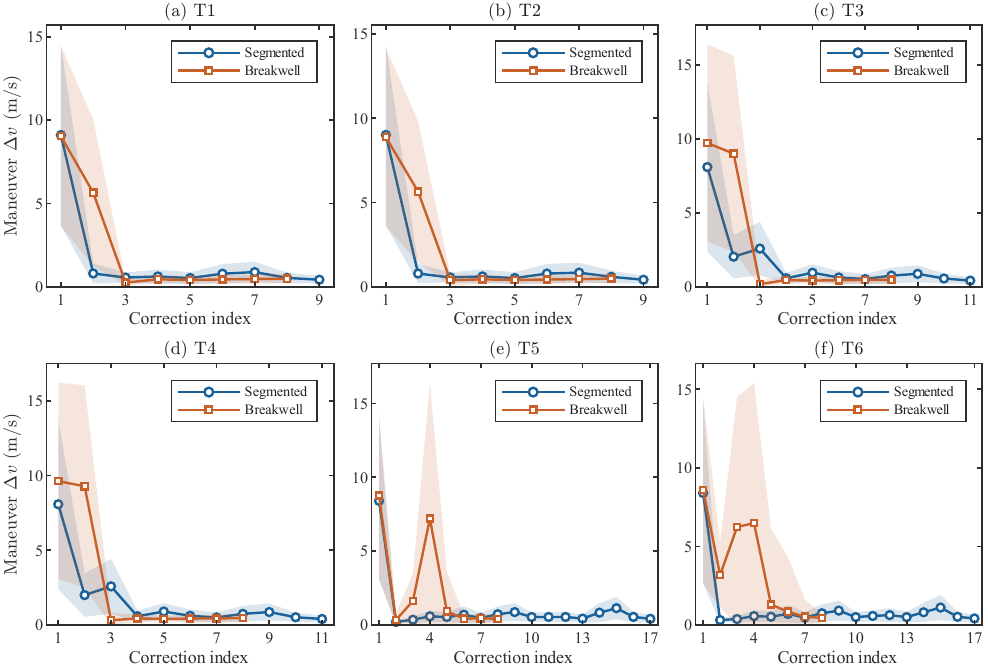}
	\caption{Per-correction maneuver statistics on common-success pairs. Lines show the mean pulse magnitude, and shaded bands show mean $\pm$ one standard deviation.}
	\label{fig:pulse_statistics}
\end{figure}

The lower and more stable total cost can be traced to a different distribution of maneuver effort along the transfer. The first Segmented command averages \SI{8.08}{\metre\per\second}--\SI{9.10}{\metre\per\second} across T1--T6, compared with \SI{8.59}{\metre\per\second}--\SI{9.72}{\metre\per\second} for Breakwell. Excluding the first command, the Segmented per-index means, averaged over the remaining correction indices for each trajectory, remain within \SI{0.60}{\metre\per\second}--\SI{0.98}{\metre\per\second}. Its largest per-index standard deviation ranges from \SI{5.37}{\metre\per\second} to \SI{5.75}{\metre\per\second}, compared with \SI{5.31}{\metre\per\second}--\SI{9.18}{\metre\per\second} for Breakwell (Fig.~\ref{fig:pulse_statistics}). Both statistics are slightly higher for Segmented on T1--T2.

This maneuver pattern is consistent with the error-reset mechanism observed on T3. Breakwell must occasionally recover a deviation after it has already grown over a long nonlinear interval, producing several large and disturbance-sensitive commands. The Segmented strategy instead distributes the response across successive sensitivity-defined intervals. Each correction reduces the residual passed to the next local targeting problem, so the subsequent maneuvers remain small. The proposed framework therefore exchanges a small number of potentially large recovery maneuvers for more frequent low-amplitude maintenance corrections. This increases the number of navigation and commanding opportunities, but reduces both the total impulse demand and its variability.

\subsection{Overall Performance Assessment}

Taken together, the six transfers reveal a consistent performance trade-off. Here, S/B denotes the Segmented and Breakwell strategies, respectively, and $\Delta p=p_{\mathrm S}-p_{\mathrm B}$ denotes the difference in mission-success probability. Reliability is evaluated over all paired Monte Carlo realizations, while correction cost, path RMSE, and terminal position error are evaluated only over the common-success pairs; the reported cost and path reductions are measured relative to Breakwell. The overall effect sizes are collected in Table~\ref{tab:effect_sizes}. The proposed strategy consistently improves mission reliability and, when both methods succeed, reduces correction cost and full-path deviation. Breakwell retains an advantage in conditional terminal-position accuracy, but this additional precision does not compensate for its lower recovery probability and substantially larger intermediate path deviation. The principal benefit of sensitivity-guided segmentation is therefore not minimum terminal miss distance, but reliable and fuel-efficient preservation of the complete low-energy transfer.

\begin{table}[htbp]
	\caption{Overall correction performance across the six transfers.}
	\label{tab:effect_sizes}
	\centering
	\footnotesize
	\setlength{\tabcolsep}{2.3pt}
	\renewcommand{\arraystretch}{1.08}
	\begin{tabular}{ccccccccc}
		\toprule
		& \multicolumn{2}{c}{Reliability} & & \multicolumn{2}{c}{Correction cost} & \multicolumn{2}{c}{Path adherence} & Terminal accuracy\\
		\cmidrule(lr){2-3}\cmidrule(lr){5-6}\cmidrule(lr){7-8}\cmidrule(lr){9-9}
		Case & Success S/B & $\Delta p$ & Common $n$ & Mean $D$ S/B & Reduction & Median $G_{\mathrm{path}}$ S/B & Reduction & Median $e_r$ S/B\\
		& (\%) & (pp) & & (m/s) & (\%) & (km) & (\%) & (km)\\
		\midrule
		T1 & 99.2/60.6 & 38.6 & 604 & 14.18/17.14 & 17.3 & 93.1/3552.9 & 97.4 & 12.64/0.22\\
		T2 & 98.4/60.6 & 37.8 & 599 & 14.07/17.10 & 17.7 & 90.6/3509.3 & 97.4 & 13.57/0.22\\
		T3 & 98.7/86.9 & 11.8 & 864 & 17.94/21.12 & 15.0 & 196.1/1784.6 & 89.0 & 12.08/0.25\\
		T4 & 98.9/85.5 & 13.4 & 852 & 17.85/21.43 & 16.7 & 194.2/1804.7 & 89.2 & 13.20/0.26\\
		T5 & 98.2/7.3 & 90.9 & 73 & 17.97/20.22 & 11.1 & 112.5/472.3 & 76.2 & 14.64/0.39\\
		T6 & 98.0/26.5 & 71.5 & 262 & 18.18/27.64 & 34.2 & 111.8/1114.1 & 90.0 & 11.52/0.37\\
		\bottomrule
	\end{tabular}
\end{table}

\begin{figure}[p]
	\centering
	\includegraphics[width=0.96\textwidth,height=0.88\textheight,keepaspectratio]{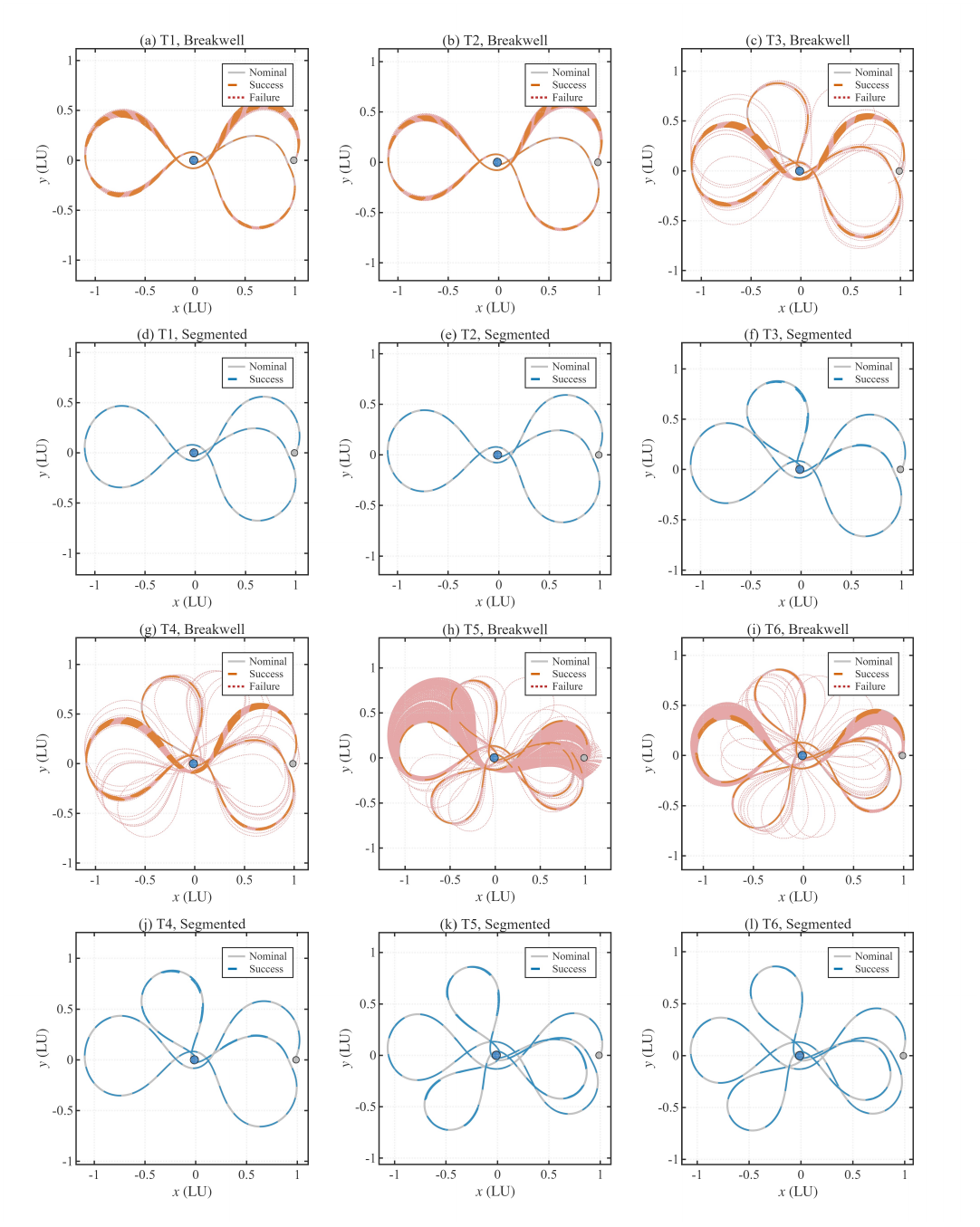}
	\caption{Paired Monte Carlo results of corrected trajectory for T1--T6.}
	\label{fig:mc_geometry_all}
\end{figure}

The Monte Carlo trajectory clouds in Figure~\ref{fig:mc_geometry_all} provide a geometric interpretation of the reliability results. The figure includes all realizations successfully recovered by the Segmented strategy and the corresponding Breakwell trajectories under identical disturbances. The Segmented trajectories generally preserve the nominal transfer geometry, while Breakwell can produce substantial departures from the intended Earth-centered loops and consequently violate the joint requirements on path recovery, terminal accuracy, or maneuver limits. For T5, 909 of the 982 Segmented-success realizations fail under Breakwell. T6 has the same segment and maneuver counts as T5 but a higher Breakwell success rate, indicating that recovery capability also depends on the timing of sensitive passages.

This behavior is consistent with the error-growth mechanism in Section~\ref{subsec5.2} and the maneuver statistics in Section~\ref{subsubsec5.3.2}. Breakwell concentrates corrections near arrival without explicitly suppressing errors before intermediate sensitive passages, allowing larger deviations to propagate downstream and increasing the subsequent correction demand. Although late corrections may improve terminal accuracy, they cannot recover deviations that have already developed along the intermediate path. By placing corrections around the sensitivity boundaries, the Segmented strategy limits this error accumulation and residual carryover. It therefore achieves higher joint mission reliability while reducing total correction $\Delta v$ and improving path preservation on common-success pairs.

\FloatBarrier

\FloatBarrier

\section{Conclusions}\label{sec6}

This paper developed an LLE-guided segmented mid-course correction framework for low-energy Earth-Moon transfers with strongly nonuniform error growth. The instantaneous LLE was used to identify sensitive passages and divide the nominal trajectory into correction segments. The LLE maxima were adopted as segment boundaries and downstream targets, thereby determining the maneuver architecture from the local dynamical structure of each transfer. A segment-wise parameterization reduced the correction schedule to one timing variable per segment, while second-order state transition tensors combined with a Levenberg--Marquardt iteration were used for nonlinear targeting. The remaining timing variables were optimized under injection, navigation, and maneuver-execution uncertainties using a chance-constrained sample-average approximation. This framework separates dynamics-informed architecture design from uncertainty-aware timing optimization and converts a complex multi-encounter correction problem into a sequence of lower-dimensional local targeting problems.

Paired Monte Carlo simulations on six representative transfers showed that the proposed strategy achieved joint success rates of 98.0--99.2\%, improving trajectory-recovery reliability by 11.8--90.9 percentage points relative to the classical Breakwell strategy. On the common-success realizations, the mean correction cost was reduced by 11.1--34.2\%, and the median full-path position RMSE was reduced by 76.2--97.4\%. Although Breakwell produced smaller terminal position errors among the common-success pairs, the proposed strategy satisfied the prescribed arrival requirements while maintaining the spacecraft substantially closer to the intended low-energy path. Its principal benefit is therefore not the minimization of terminal miss distance, but the reliable and fuel-efficient recovery of the complete transfer trajectory. The results further indicate that correction difficulty depends on the location and sequence of sensitive passages rather than on flight time or maneuver count alone. Future work will extend the framework to three-dimensional high-fidelity ephemeris dynamics and more realistic uncertainty models. Adaptive maneuver management, including maneuver merging and skipping, will also be investigated.

\printcredits

\section*{Declaration of competing interest}
The authors declare that they have no known competing financial interests or personal relationships that could have appeared to influence the work reported in this paper.

\section*{Acknowledgements}

This work was supported by the National Natural Science Foundation of China (Grant No. 12525204), the Beihang University Fund Program for Talent Development Support Scheme (Grant No. 4263023), and the State Key Laboratory Project of High-Efficiency Reusable Aerospace Transportation Technology (Grant No. RATT-ZZ03-2025029). The authors thank Dai Shen of Beihang University for helpful comments.

\bibliographystyle{elsarticle-num}
\bibliography{cas-refs}


\end{document}